\documentclass[a4paper,11pt]{amsart}
\usepackage{amssymb}
\usepackage{graphicx}
\usepackage{rotating}

\def\bfit{\bfseries\itshape}

\input xypic
\xyoption{all}
\xyoption{arc}

\newtheorem{theo}{Theorem}[section]

\newtheorem{prop}[theo]{Proposition}

\newtheorem{lem}[theo]{Lemma}

\newtheorem{coro}[theo]{Corollary}

\def\remark#1{{\refstepcounter{theo}\label{#1}\noindent\sc Remark  
\arabic{section}.\arabic{theo} - }}
\def\example#1{{\refstepcounter{theo}\label{#1}\noindent\sc Example 
\arabic{section}.\arabic{theo} - }}

\def\equat{\refstepcounter{theo}$$~}
\def\endequat{\leqno{\boldsymbol{(\arabic{section}.\arabic{theo})}}~$$}

\newcounter{soussection}[section]

\def\soussection#1{\refstepcounter{soussection}
\noindent{\bfit \arabic{section}.\Alph{soussection}. #1.}}

    \def\CM{{\mathbb{C}}}

    \def\FM{{\mathbb{F}}}

    \def\NM{{\mathbb{N}}}

    \def\QM{{\mathbb{Q}}}
    \def\RM{{\mathbb{R}}}
\def\SG{{\mathfrak S}}    
\def\TG{{\mathfrak T}}

    \def\ZM{{\mathbb{Z}}}

\def\Cb{{\mathbf C}}    \def\CC{{\mathcal{C}}}
    \def\DC{{\mathcal{D}}}
    \def\EC{{\mathcal{E}}}
    \def\FC{{\mathcal{F}}}

\def\Ib{{\mathbf I}}    \def\IC{{\mathcal{I}}}
    
    \def\KC{{\mathcal{K}}}

    \def\NC{{\mathcal{N}}}
    
    \def\PC{{\mathcal{P}}}

    \def\TC{{\mathcal{T}}}

    \def\WC{{\mathcal{W}}}
    \def\XC{{\mathcal{X}}}

  \def\xti{{\tilde{x}}}

          \def\aba{{\bar{a}}}

\def\a{\alpha}

\def\g{\gamma}
\def\G{\Gamma}
\def\d{\delta}
\def\D{\Delta}
\def\e{\varepsilon}
\def\ph{\varphi}
\def\l{\lambda}
\def\L{\Lambda}

\def\s{\sigma}
\def\th{\theta}

\def\t{\tau}
\def\z{\zeta}

        \def\Delt{{\tilde{\Delta}}}

\def\lamb{{\boldsymbol{\lambda}}}       
\def\Lamb{{\boldsymbol{\Lambda}}}

          \def\lamh{{\hat{\lambda}}}
          
               \def\muh{{\hat{\mu}}}

\DeclareMathOperator{\End}{{\mathrm{End}}}

\DeclareMathOperator{\Hom}{{\mathrm{Hom}}}

\DeclareMathOperator{\Id}{{\mathrm{Id}}}

\DeclareMathOperator{\im}{{\mathrm{Im}}}

\DeclareMathOperator{\Ind}{{\mathrm{Ind}}}

\DeclareMathOperator{\Irr}{{\mathrm{Irr}}}
\DeclareMathOperator{\Ker}{{\mathrm{Ker}}}

\DeclareMathOperator{\Rad}{{\mathrm{Rad}}}
\DeclareMathOperator{\reg}{{\mathrm{reg}}}

\DeclareMathOperator{\Res}{{\mathrm{Res}}}

\DeclareMathOperator{\Supp}{{\mathrm{Supp}}}

\DeclareMathOperator{\rank}{{\mathrm{rk}}}
\DeclareMathOperator{\Cartan}{{\mathrm{Cartan}}}

\def\to{\rightarrow}
\def\longto{\longrightarrow}
\def\injto{\hookrightarrow}

\def\fonction#1#2#3#4#5{\begin{array}{rccc}
{#1} : & {#2} & \longto & {#3} \\
& {#4} & \longmapsto & {#5} 
\end{array}}

\def\vide{\varnothing}

\def\DS{\displaystyle}
\def\SS{\scriptstyle}

\def\finl{~$\SS \square$}

\def\infspe{\hspace{0.1em}\mathop{\preccurlyeq}\nolimits\hspace{0.1em}}

\def\lexp#1#2{\kern\scriptspace\vphantom{#2}^{#1}\kern-\scriptspace#2}
\def\le{\hspace{0.1em}\mathop{\leqslant}\nolimits\hspace{0.1em}}
\def\ge{\hspace{0.1em}\mathop{\geqslant}\nolimits\hspace{0.1em}}

\mathchardef\inferieur="321E
\mathchardef\superieur="321F

\def\eqna{\begin{eqnarray*}}
\def\endeqna{\end{eqnarray*}}

\def\rem{\noindent{\sc{Remark}~-} }

\def\para{parabolic subgroup }

\def\itemth#1{\item[${\mathrm{(#1)}}$]}

\def\Cartan{{\mathrm{Cartan}}}

\def\VDash{{\, \scriptstyle{ |\models}\,}}
\DeclareMathOperator{\Comp}{{\mathrm{Comp}}}
\DeclareMathOperator{\Bip}{{\mathrm{Bip}}}
\DeclareMathOperator{\cox}{{\mathrm{cox}}}
\def\para{{\mathrm{par}}}

\DeclareMathOperator{\Sat}{{\mathrm{Sat}}}
\DeclareMathOperator{\val}{{\mathrm{val}}}
\DeclareMathOperator{\length}{{\mathrm{lg}}}
\DeclareMathOperator{\Mat}{{\mathrm{Mat}}}
\def\ubar{{\bar{1}}}
\def\dbar{{\bar{2}}}
\def\tbar{{\bar{3}}}

\begin{document}

\baselineskip=16pt

\title{Representation theory of \\ Mantaci-Reutenauer algebras}

\author{C. Bonnaf\'e}
\address{\noindent 
Labo. de Math. de Besan\c{c}on (CNRS: UMR 6623), 
Universit\'e de Franche-Comt\'e, 16 Route de Gray, 25030 Besan\c{c}on
Cedex, France} 

\makeatletter
\email{bonnafe@math.univ-fcomte.fr}

\makeatother


\date{\today}

\begin{abstract} 
We study some aspects of the representation theory 
of Mantaci-Reutenauer algebras: Cartan matrix, 
Loewy length, modular representations. 
\end{abstract}

\maketitle

\pagestyle{myheadings}

\markboth{\sc C. Bonnaf\'e}{\sc Mantaci-Reutenauer algebra}

Let $W_n$ be a Coxeter group of type $B_n$ (i.e. the group of permutations 
$\s$ of $I_n=\{\pm 1,\dots,\pm n\}$ such that $\s(-i)=-\s(i)$ for every $i \in I_n$) 
and let $R$ be a commutative ring. Mantaci and Reutenauer \cite{mantaci} have defined 
a subalgebra $R\Sigma'(W_n)$ of the group algebra $RW_n$ which contains 
both the Solomon descent algebra of the symmetric group $\SG_n$ and the one 
of $W_n$. In \cite{BH}, the authors have provided another construction 
of the Mantaci-Reutenauer algebra $R\Sigma'(W_n)$ which relies more 
on the structure of $W_n$ as a Coxeter group. As a consequence of their work, 
they were able to generalize to this algebra the classical results 
of Solomon on the Solomon descent algebra (construction of a morphism 
to the character ring of $W_n$, description of the radical whenever 
$R$ is a field of characteristic $0$...). For instance, the description 
of the simple $\QM\Sigma'(W_n)$-modules was obtained in \cite[Proposition 3.11]{BH}: 
they are all of dimension $1$. 

In this paper, we study the representation theory of $K\Sigma'(W_n)$ 
whenever $R=K$ is a field of any characteristic: simple modules, radical, 
projective modules, Cartan matrix... We also define some 
morphisms between different Mantaci-Reutenauer algebras. 
Let us gather here some of the main results obtained all along the text:

\bigskip

\noindent{\bf Theorem.} 
{\it Let $p$ denote the characteristic of $K$. Then:
\begin{itemize}
\itemth{a} There exists a natural morphism of algebras $K\Sigma'(W_n) \to K\Sigma'(W_{n-1})$; 
it is surjective if $p=0$.

\itemth{b} If $p \neq 2$, then the Loewy length of $K\Sigma'(W_n)$ is $n$. 
If $p=2$, then this Loewy length lies in $\{n,n+1,\dots,2n-1\}$. 

\itemth{c} If $p$ does not divide $|W_n|$ (i.e. is $p=0$ or $p > \max(2,n)$), 
then the Cartan matrix of $K\Sigma'(W_n)$ is unitriangular.

\itemth{d} If $p$ does not divide the order of $W_n$, 
then the Cartan matrix of $K\Sigma'(W_n)$ 
is a submatrix of the Cartan matrix of $K\Sigma'(W_{n+1})$. 
\end{itemize}}

\bigskip

Note that, in the statement (a), we expect that the homomorphism 
is surjective even if $p > 0$, but we are unable to prove it. 
In the statement (b), we expect that the Loewy length of $\FM_{\! 2} \Sigma'(W_n)$ 
is equal to $2n-1$ whenever $n \ge 2$. 

\medskip

The paper is organized as follows. In the first section, we gather 
(and sometimes improve, or make more precise) 
some of the principal results of \cite{BH}. In the second 
section, we study some particular families of left, 
right and two-sided ideals of $R\Sigma'(W_n)$. In the third 
section, we introduce a class of {\it positive} elements of 
$K\Sigma'(W_n)$ (whenever $K$ is an ordered field) and study 
the ideals they generate (and also some other properties: 
centralizer, minimal polynomial). In the fourth section, 
we study the action of the longest element $w_n$ of $W_n$ on 
simple modules and on $K\Sigma'(W_n)$: since $w_n$ is central 
(and is an element of $K\Sigma'(W_n)$), this provides 
a first decomposition of the Mantaci-Reutenauer algebra 
(at least when $K$ is not of characteristic $2$: we also 
give a basis of $K\Sigma'(W_n)$ consisting of eigenvectors 
for the action of $w_n$ by left multiplication). In the fifth 
section, we define some morphisms between Mantaci-Reutenauer 
algebras and prove the statement (a) of the Theorem above. 
In the sixth section, we study the simple modules 
and compute explicitly the radical of $K\Sigma'(W_n)$ (this 
is done in any characteristic). Section 7 is devoted to 
the computation of the Loewy length of $K\Sigma'(W_n)$, 
that is to the proof of the statement (b) of the above Theorem. 
We also obtain the Loewy length of the algebra $K\Irr W_n$ 
in any characteristic. The section 8 is concerned with 
the projective modules and the Cartan matrix of $K\Sigma'(W_n)$: 
the statement (c) and (d) of the above Theorem are proved. 
We also obtain some results about the structure of $K W_n$ 
as a $K\Sigma'(W_n)$-module. We give in section 9 
some numerical results (character tables, primitive idempotents 
and the Cartan matrices for small values of $n$). 
In the final section, we address some questions that are raised by 
the present work. 

Most of this work is largely inspired by works of several 
authors on Solomon descent algebras (see for instance \cite{atkinson}, 
\cite{bbht}, \cite{apv}, \cite{BP}...). Sections 2, 3, 4, 5 
are analogous to \cite[\S 2, 3, 4]{BP} (for \S 5, see also \cite{atkinson} 
and \cite{bbht}). Sections 6 and 8 are inspired by \cite{apv}). 
Section 7 is the analogue of \cite[\S 5]{BP}. The question 6 in section 10 
has been suggested by a similar question of G. Pfeiffer on Solomon 
descent algebras.

\tableofcontents

\section{Notation, preliminaries}

\medskip

\soussection{General notation}
All along this paper, $R$ will denote a fixed commutative 
ring and $K$ a fixed field. If $G$ is a finite group, the group 
algebra of $G$ over $R$ is denoted by $RG$ and the set of 
irreducible characters of $G$ over $\CM$ is denoted by $\Irr G$. 
We denote by $R\Irr G$ the ring of formal $R$-linear combinations 
of irreducible characters of $G$ (with multiplication given by 
tensor product). In particular, $\ZM\Irr G$ can 
be identified with the Grothendieck ring of the category of finite 
dimensional $\CM G$-module (which is usually called the {\it character 
ring} of $G$) and $R\Irr G=R \otimes_\ZM \ZM\Irr G$. 
If $A$ is a finite dimensional $K$-algebra, we denote by $\Rad A$ 
its radical. 

\bigskip

\soussection{Weyl group of type ${\boldsymbol{B_n}}$\label{section type B}}
If $n \ge 1$, we denote by $(W_n,S_n)$ a Weyl group of type $B_n$: 
write $S_n=\{t,s_1,s_2,\dots,s_{n-1}\}$ in such a way that the Dynkin 
diagram of $W_n$ is
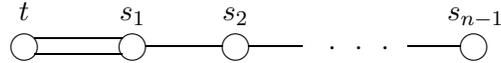
\begin{center}
\begin{picture}(220,30)
\put( 40, 10){\circle{10}}
\put( 44,  7){\line(1,0){33}}
\put( 44, 13){\line(1,0){33}}
\put( 81, 10){\circle{10}}
\put( 86, 10){\line(1,0){29}}
\put(120, 10){\circle{10}}
\put(125, 10){\line(1,0){20}}
\put(155,  7){$\cdot$}
\put(165,  7){$\cdot$}
\put(175,  7){$\cdot$}
\put(185, 10){\line(1,0){20}}
\put(210, 10){\circle{10}}
\put( 38, 20){$t$}
\put( 76, 20){$s_1$}
\put(115, 20){$s_2$}
\put(200, 20){$s_{n{-}1}$}
\end{picture}
\end{center}
Let $S_{-n}=\{s_1,s_2,\dots,s_{n-1}\}$ and $W_{-n}=<S_{-n}>$. Note that 
$W_{-n} \simeq \SG_n$. We denote by $\ell : W_n \to \NM$ the 
length function attached to $S_n$. 

Let $I_n=\{\pm 1,\pm 2,\dots,\pm n\}$. We identify $W_n$ with the 
group of permutations $\s$ of $I_n$ such that $\s(-i)=-\s(i)$ 
for every $i \in I_n$. The identification is as follows: 
$t$ corresponds to the transposition $(1,-1)$ while $s_i$ 
corresponds to $(i,i+1)(-i,-i-1)$. 
Let $t_1=t$ and, if $1 \le i \le n-1$, let $t_{i+1}=s_it_is_i$. 
As a permutation of $I_n$, $t_i$ is equal to $(i,-i)$. 
Now, we set $T_n=\{t_1,\dots,t_n\}$ and $S_n'=S_n \cup T_n$. Then the
reflection subgroup $\TG_n$ generated by $T_n$ is naturally
identified with $(\ZM/2\ZM)^n$. Therefore $W_n=W_{-n} \ltimes
\TG_n$ is, abstractly, the wreath product of $\SG_n$ by  $\ZM/2\ZM$. 

Let $(e_1,\dots,e_n)$ denote the canonical basis of the euclidean 
$\RM$-vector space $\RM^n$. If $\a \in \RM^n$, we denote by $s_\a$ the 
orthogonal reflection such that $s_\a(\a)=-\a$. Let 
$$\Phi_n=\{\pm 2e_i~|~1 \le i \le n\} \cup \{ \pm e_i \pm e_j~|~1 \le i < j \le n\}.$$
Then $\Phi_n$ is a root system and $W_n$ can be also identified 
with the Weyl group of $\Phi_n$: through this identification, we 
have $t_i=s_{2e_i}$ and $s_i = s_{e_{i+1}-e_i}$. Let 
$$\D_n=\{2e_1,e_2-e_1,e_3-e_2,\dots,e_n-e_{n-1}\}.$$
Then $\D_n$ is a basis of $\Phi_n$ and we denote by $\Phi_n^+$ the 
set of roots which are linear combinations with non-negative coefficients 
of roots in $\D_n$. If $\a \in \Phi_n$, we write $\a > 0$ if $\a \in \Phi_n^+$ 
and $\a < 0$ otherwise. 

\bigskip

\soussection{Signed compositions, bipartitions\label{section signed}} 
A {\it signed composition} is a finite 
sequence $C = (c_1, \dots, c_r)$  of non-zero elements of $\ZM$. The number $r$ is
called the {\it length} of $C$ and will be denoted by $\length(C)$. 
We denote by $\length^+(C)$ (respectively $\length^-(C)$) the number 
of positive (respectively negative) parts of $C$. In particular, 
$\length(C)=\length^+(C)+\length^-(C)$. 
We set $|C| =\sum_{i=1}^r | c_i |$. If $|C|=n$,
we say that $C$ is a {\it signed composition of $n$} and we write $C \VDash n$.
We also define
$C^+=(|c_1|,\dots,|c_r|) \VDash n$ and $C^-=-C^+$. We denote
by $\Comp(n)$ the set of signed compositions of $n$. In
particular, any composition is a signed composition (any part is
positive). Note that 
\equat\label{nombre composition}
|\Comp(n)|=2.3^{n-1}.
\endequat
If $C=(c_1,\dots,c_r)$ and $D=(d_1,\dots,d_s)$ are signed compositions of 
$m$ and $n$ respectively, we denote by $C \sqcup D$ the signed 
composition $(c_1,\dots,c_r,d_1,\dots,d_s)$ of $m+n$.

A \textit{bipartition} of $n$ is a pair
$\l= (\l^+ , \l^-)$ of partitions such that $|\l|:=|\l^+|+|\l^-|=n$. 
We set $\length^+(\l)=\length(\l^+)$, $\length^-(\l)=\length(\l^-)$ and 
$\length(\l)=\length(\l^+) + \length(\l^-)$. 
We write $\l \Vdash n$ to say that $\l$ is a bipartition of $n$, and
the set of bipartitions of $n$ is denoted by
$\Bip(n)$. 
We define ${\hat{\l}}$ as the signed composition of $n$ obtained by
concatenation of $\l^+$ and $-\l^-$. In other words, 
$\lamh=\l^+ \sqcup -\l^-$. The map 
$\Bip(n) \to \Comp(n)$, $\l \mapsto \hat{\l}$ is injective.

Now, let $C$ be a signed composition of $n$. We define
$\lamb(C)=(\l^+,\l^-)$ as the bipartition of $n$ such that $\l^+$
(resp. $\l^-$) is obtained from $C$ by reordering if necessary 
the positive parts of $C$ (resp. the absolute value of the
negative parts of $C$). Note that $\length(\lamb(C))=\length(C)$, 
$\length^+(\lamb(C))=\length^+(C)$ and $\length^-(\lamb(C))=\length^-(C)$. 
One can easily check that the map
$$\lamb : \Comp(n) \longto \Bip(n)$$
is surjective (indeed, if $\l \in \Bip(n)$, then
$\lamb(\hat{\l})=\l$).

\bigskip 

\soussection{A class of reflection subgroups of ${\boldsymbol{W_n}}$\label{class}}
Now, to each $C=(c_1, \dots, c_r) \VDash n$, we associate a reflection 
subgroup $W_C$ of $W_n$ which is isomorphic to 
$W_{c_1} \times \ldots \times W_{c_r}$. We proceed 
as follows:  for $1\leq i\leq r$, set
$$I_C^{(i)}=\begin{cases}
            I_{C,+}^{(i)} & \text{if } c_i < 0,\\
        I_{C,+}^{(i)} \cup -I_{C,+}^{(i)} & \text{if } c_i > 0,
        \end{cases}$$
where $I_{C,+}^{(i)}$ is the set of natural numbers $k$ such that 
$|c_1|+\dots+|c_{i-1}|+1 \le k \le |c_1|+\dots+|c_i|$. Then
$$W_C=\{w \in W_n ~|~ \forall~1 \le i \le r,~w(I_C^{(i)})=I_C^{(i)}\}$$
is a reflection subgroup generated by
$$S_C=  (S_{-n} \cap W_C) \cup \{t_{|c_1| + \dots+|c_{j-1}|+1} \in T_n ~|~
c_j >0\} \qquad \subset S_n'.$$
Note that $(W_C,S_C)$ is a Coxeter group. 
Moreover, $W_C \simeq W_{c_1} \times \dots \times W_{c_r}$. 
Let $T_C=T_n \cap W_C$. Then $W_C=\SG_{C^+} \ltimes <T_C>$, where 
$\SG_{C^+}=W_{C^-}$.

Now let $\Phi_C=\{\a \in \Phi_n~|~s_\a \in W_C\}$. Then $\Phi_C$ is a root 
system and $W_C$ is naturally identified with the Weyl group of $\Phi_C$. 
Let $\Phi_C^+=\Phi_C \cap \Phi_n^+$ and $\D_C=\{\a \in \Phi_C^+~|~s_\a \in S_C\}$. 
Then $\D_C$ is a basis of $\Phi_C$, and $\Phi_C^+$ 
is a positive root system of $\Phi_C$.


If $C, D\VDash n$, then we write $C \subset D$ if $W_C \subset
W_D$. This defines an order $\subset$ on $\Comp(n)$. 

\bigskip

\remark{proprietes inclusion} If $C \subset D$, 
then $\length(C) \ge \length(D)$ and $\length^-(C) \ge \length^-(D)$. 
If $C \subset D$, $\length(C) =\length(D)$ and $\length^-(C)=\length^-(D)$, 
then $C=D$.\finl

\bigskip

\noindent{\sc Example - } It might happen that $C \subset D$ and 
$\length^+(C) < \length^+(D)$. For example, take $C=(-n)$ and $D=(n)$.\finl

\bigskip

\soussection{Conjugacy classes\label{section conjugacy}}
If $C \VDash n$, we denote by $\cox_C$ a Coxeter element of $(W_C,S_C)$. 
If $C$, $C' \subset D$ and if $W_C$ and $W_{C'}$ are
conjugate under $W_D$, then we write $C \equiv_D C'$. Note that 
$\cox_C$ and $\cox_{C'}$ are conjugate in $W_D$ if and only if $C \equiv_D C'$. 
Moreover, every element of $W_D$ is $W_D$-conjugate to $\cox_C$ for 
some $C \subset D$. If $D=(n)$, we write 
$\equiv$ instead of $\equiv_D$. We recall the following easy proposition:

\bigskip

\noindent{\bf Proposition A.}
{\it Let $C,D\VDash n$. Then $W_C$ and $W_D$ are conjugate in $W_n$ if
and only if $\lamb(C)=\lamb(D)$.}

\bigskip

If $w \in W_n$, we denote by $\Lamb(w)$ the unique bipartition $\l$ 
of $n$ such that $w$ is conjugate to $\cox_C$ for some 
(every) $C \in \lamb^{-1}(\l)$. The map 
$$\Lamb : W_n \longto \Bip(n)$$
is well-defined, surjective and its fibers are precisely the 
conjugacy classes of $W_n$: if $\l \in \Bip(n)$, we set 
$\CC(\l)=\Lamb^{-1}(\l)$ and we fix an element $\cox_\l \in \CC(\l)$ 
(if $C \in \Comp(n)$, $\cox_{\lamb(C)}$ is conjugate to $\cox_C$). 
We denote by $o(\l)$ the order of an element 
of $\CC(\l)$: if $\l=(\l^+,\l^-)$ where 
$\l^+=(\l_1^+,\dots,\l_k^+)$ and $\l^-=(\l_1^-,\dots,\l_l^-)$, 
then $o(\l)$ is the least common multiple of 
$\{2\l_1^+,\dots,2\l_k^+,\l_1^-,\dots,\l_l^-\}$. 

\bigskip

\soussection{Mantaci-Reutenauer algebra}
Let $C \VDash n$, then
$$
X_C=\{x \in W_n~|~\forall~w \in W_C,~\ell(xw) \ge \ell(x)\}
$$
is a distinguished set of {\it minimal coset representatives}
for $W_n/W_C$ (see \cite[Proposition 2.8 (a)]{BH}). It is readily seen that 
\eqna
X_C&=& \{w \in W_n~|~\forall~s \in S_C,~\ell(ws)>\ell(w)\}\\
&=& \{w \in W_n~|~\forall~\a \in \Phi_C^+,~w(\a) > 0\} \\
&=& \{w \in W_n~|~\forall~\a \in \D_C,~w(\a) > 0\}.
\endeqna
Now, we set
$$x_C=\sum_{w \in X_C} w \qquad \in R W_n.$$
(Recall that $R$ is a fixed commutative ring.) 
By \cite[\S 3.1]{BH}, the family $(x_C)_{C \in \Comp(n)}$ 
is free over $R$. Let
$$R\Sigma'(W_n)=\mathop{\oplus}_{C \in \Comp(n)} R x_C \qquad
\subset RW_n.$$
For simplification, we set $\Sigma'(W_n)=\ZM \Sigma'(W_n)$, so 
that $R\Sigma'(W_n)=R \otimes_\ZM \Sigma'(W_n)$. 

\bigskip

\rem The algebra $R\Sigma'(W_n)$ is nothing else but the algebra constructed by Mantaci 
and Reutenauer \cite{mantaci} by combinatorial methods 
(see \cite[Remark of Subsection 3.1]{BH} 
for the identification).\finl

\bigskip

Let $(\xi_C)_{C \in \Comp(n)}$ denote the basis $\Hom_R(R\Sigma'(W_n),R)$ dual 
to $(x_C)_{C \in \Comp(n)}$. 
In other words, we have, for every $x \in R\Sigma'(W_n)$, 
$$x=\sum_{C \in \Comp(n)} \xi_C(x) x_C.$$
We now define
$$\th_n^R : R\Sigma'(W_n) \longto R \Irr W_n$$
as the unique $R$-linear map such that 
$$\th_n^R(x_C)=\Ind_{W_C}^{W_n} 1_C$$
for every $C \in \Comp(n)$.
Here, $1_C$ is the trivial character of $W_C$. We denote by $\e_C$
the sign character of $W_C$.
We can now recall the following result.

\bigskip

\noindent{\bf Theorem B [BH, Theorem 3.7].} 

{\it 
\begin{itemize}
\itemth{a} $R\Sigma'(W_n)$ is a unitary sub-$R$-algebra of $R W_n$. 

\itemth{b} $\th_n^R : R\Sigma'(W_n) \to R \Irr W_n$ is a 
morphism of $R$-algebras.

\itemth{c} $\th_n^R$ is surjective and 
$\Ker \th_n^R=\DS{\sum_{C \equiv D} R(x_C - x_D)}$. 

\itemth{d} If $K$ is a field of characteristic $0$, then 
$\Ker \th_n^K$ is the radical of the $K$-algebra $K\Sigma'(W_n)$. 
\end{itemize}}

\bigskip

Let $\Comp^+(n)$ be the set of compositions of $n$. A signed 
composition $C=(c_1,\dots,c_r)$ is called {\it semi-positive} 
(resp. {\it parabolic}) if $c_i \ge -1$ (resp. $c_i < 0$) for every $i \ge 1$ 
(resp. for every $i \ge 2$). Note that $C$ is parabolic if and only if $W_C$ is a 
standard parabolic subgroup of $W$ (i.e. if and only if $S_C \subset S_n$). 
We denote by $\Comp_\para(n)$ the set of parabolic compositions of 
$n$. Let
$$R\Sigma(W_n)=\mathop{\oplus}_{C \in \Comp_\para(n)} R x_C$$
$$R\Sigma(\SG_n) = \mathop{\oplus}_{C \in \Comp^+(n)} R x_C.\leqno{\text{and}}$$
Then $R\Sigma(W_n)$ and $R\Sigma(\SG_n)$ are sub-$R$-algebras of $R\Sigma'(W_n)$: 
$R\Sigma(W_n)$ is the Solomon descent algebra of $W_n$ (see \cite{solomon} 
for the definition of Solomon descent algebras of finite Coxeter groups) 
and it is easy to check \cite[\S 3.2]{BH} that 
$R\Sigma(\SG_n)$ is the Solomon descent algebra of $\SG_n=W_{-n}$. 

The restriction of $\th_n^R$ to $R\Sigma(W_n)$ is equal to 
the classical Solomon homomorphism. On the other hand, the canonical surjective 
morphism $W_n \to \SG_n$ induces an injective morphism of algebras 
$R\Irr \SG_n \injto R \Irr W_n$. We view $R\Irr \SG_n$ naturally as a subalgebra 
of $R\Irr W_n$ through this morphism. Then the image, through $\th_n^R$, of an element 
of $R\Sigma(\SG_n)$ belongs to $R\Irr \SG_n$. Also, the restriction 
of $\th_n^R$ to a morphism (still denoted by $\th_n^R$) of algebras 
$R\Sigma(\SG_n) \to R\Irr \SG_n$ is again equal to the classical Solomon 
homomorphism. By construction, the diagram
\equat\label{diagramme sn}
\diagram
R\Sigma(\SG_n) \xto[0,2]|<\ahook
\ddto_{\DS{\th_n^R}} && R\Sigma'(W_n) \ddto^{\DS{\th_n^R}} \\
&& \\
R\Irr \SG_n \xto[0,2]|<\ahook
&& R \Irr W_n 
\enddiagram
\endequat
is commutative \cite[Diagram 3.4]{BH}. 

\bigskip

\soussection{On the multiplication in ${\boldsymbol{R\Sigma'(W_n)}}$}
By Theorem B, $R\Sigma'(W_n)$ is a sub-$R$-algebra of $R W_n$ 
and $\th_n^R$ is a morphism of algebras. However, the multiplication 
in $R\Sigma'(W_n)$ is not described. In fact, it turns out 
that its description is much more complicated than the multiplication 
in the Solomon descent algebra. Theoretically, it is possible to 
extract from the proof of \cite[Theorem 3.7]{BH} an inductive 
process for this multiplication. We shall not do it here. 
We shall just give some easy consequences of this inductive process. 

First, if $\FC$ is a subset of $\Comp(n)$, we set
$$R\Sigma_\FC'(W_n) = \mathop{\oplus}_{C \in \FC} R x_C.$$
For instance, $R\Sigma(W_n)=R\Sigma_{\Comp_\para(n)}'(W_n)$ and 
$R\Sigma(\SG_n)=R\Sigma_{\Comp^+(n)}'(W_n)$. 

We shall now describe an order $\infspe$ on $\Comp(n)$ which is finer than 
$\subset$. Let $C$ and $D$ 
be two signed composition of $n$. We write $C \infspe D$ if 
one of the following two conditions is satisfied:

\begin{quotation}
\noindent (1) $C \subset D$.

\noindent (2) $C \subset D^+$ and $\length(C) > \length(D)$ and 
$\length^-(C) \ge \length^-(D)$.
\end{quotation}

One can easily check that it defines an order $\infspe$ on $\Comp(n)$ 
(see Remark \ref{proprietes inclusion}). We set 
$$\FC_{\inferieur D}=\{C \in \Comp(n)~|~C \inferieur D\}$$
$$\FC_{\infspe D}=\{C \in \Comp(n)~|~C \infspe D\}.\leqno{\text{and}}$$
For simplification, we set
$$R\Sigma_{\inferieur D}'(W_n)=R\Sigma_{\FC_{\inferieur D}}'(W_n)
\text{ and }R\Sigma_{\infspe D}'(W_n)=R\Sigma_{\FC_{\infspe D}}'(W_n).$$

\bigskip

\remark{comparaison} 
If $C \infspe D$, then $C^+ \subset D^+$, $\length(C) \ge \length(D)$ 
and $\length^-(C) \ge \length^-(D)$.\finl

\bigskip 

We shall now describe a preorder $\subset_\lamb$ on $\Comp(n)$. 
First, note that the order $\subset$ on $\Comp(n)$ induces an order on 
$\Bip(n)$ which we still denote by $\subset$. If $C$ and 
$D$ are two signed compositions of $n$, we then write $C \subset_\lamb D$ 
if $\lamb(C) \subset \lamb(D)$. In other words, $C \subset_\lamb D$ 
if and only if $W_C$ is contained in some conjugate of $W_D$. 
We write $C \varsubsetneq_\lamb D$ if $\lamb(C) \varsubsetneq \lamb(D)$. 
Similarly as above, we set
$$\FC_{\varsubsetneq_\lamb D}=\{C \in \Comp(n)~|~C \varsubsetneq_\lamb D\}$$
$$\FC_{\subset_\lamb D}=\{C \in \Comp(n)~|~C \subset_\lamb D\}.\leqno{\text{and}}$$
For simplification, we set
$$R\Sigma_{\varsubsetneq_\lamb D}'(W_n)=R\Sigma_{\FC_{\varsubsetneq_\lamb D}}'(W_n)
\qquad\text{and}\qquad
R\Sigma_{\subset_\lamb D}'(W_n)=R\Sigma_{\FC_{\subset_\lamb D}}'(W_n)$$

\bigskip

\remark{preordre}
It is easily checked that $\subset_\lamb$ is a preorder on $\Comp(n)$ and 
that the equivalence relation associated to the preorder 
$\subset_\lamb$ is exactly the relation $\equiv$.\finl

\bigskip

We now recall some notation from \cite[Proposition 2.13]{BH}. 
If $C$ and $D$ are two signed compositions of $n$, we set 
$$X_{CD}=X_C^{-1} \cap X_D.$$
Moreover, if $d \in X_{CD}$, we denote by $C \cap \lexp{d}{D}$ the unique 
signed composition of $n$ such that $W_C \cap \lexp{d}{W_D}=W_{C \cap \lexp{d}{D}}$. 
If $C$, $C' \subset D$, we set $X_C^D=X_C \cap W_D$, 
$x_C^D=\sum_{w \in X_C^D} w \in RW_D$, $X_{CC'}^D=X_{CC'} \cap W_D$, 
$R\Sigma'(W_D)=\mathop{\oplus}_{C \subset D} R x_C^D$ and we define 
$\th_D^R : R\Sigma'(W_D) \to R\Irr W_D$, $x_C^D \mapsto \Ind_{W_C}^{W_D} 1_C$. 
Then $R\Sigma'(W_D)$ is a sub-$R$-algebra of $RW_D$ and $\th_D^R$ is a surjective 
morphism of algebras. Moreover, if $D=(d_1,\dots,d_r)$, then 
$$R\Sigma'(W_D) \simeq R\Sigma'(W_{d_1}) \otimes_R \dots \otimes_R 
R\Sigma'(W_{d_r}),$$
where $\Sigma'(W_d)=\Sigma(\SG_{-d})$ if $d < 0$.

\bigskip

\noindent{\bf Proposition C (see [BH, Proof of Theorem 3.7]).} 
{\it Let $C$ and $D$ be two signed compositions of $n$. Then 
\begin{itemize}
\itemth{a} There is a map $f_{CD} : X_{CD} \to \Comp(n)$ such that: 
\begin{itemize}
\itemth{1} For every $d \in X_{CD}$, $f_{CD}(d) \subset D$ and 
$f_{CD}(d) \equiv_D \lexp{d^{-1}}{C} \cap D$.

\itemth{2} $x_Cx_D - \DS{\sum_{d \in X_{CD}}} x_{f_{CD}(d)} 
\in R\Sigma_{\varsubsetneq_\lamb C}'(W_n) \cap 
R\Sigma_{\inferieur D}'(W_n) \cap \Ker \th_n^R$. 
\end{itemize}

\itemth{b} If $C$ is parabolic or if $D$ is semi-positive, then 
$f_{CD}(d)=\lexp{d^{-1}}{C} \cap D$ for every $d \in X_{CD}$ and 
$x_Cx_D = \DS{\sum_{d \in X_{CD}}} x_{\lexp{d^{-1}}{C} \cap D}$. 
\end{itemize}}

\begin{proof}
In fact, (b) is proved in \cite[Example 3.2]{BH}. Let us now prove (a). 
We first need an easy lemma about double cosets representatives:
\begin{quotation}
\begin{lem}\label{double}
Let $C$, $D$ and $D'$ be three signed compositions of $n$ such that 
$D \subset D'$. Let $\EC=\{(d,e)~|~d \in X_{CD'}$ and 
$e \in X_{(\lexp{d^{-1}}{C} \cap D'),D}^{D'}\}$. Let $f : \EC \to X_{CD}$ 
be the map defined by $de \in W_C f(d,e) W_D$. Then $f$ is bijective and 
$W_{\lexp{f(d,e)^{-1}}{C} \cap D}$ is conjugate, inside $W_D$, to 
$\lexp{(de)^{-1}}{W_C} \cap W_D$. 
\end{lem}
\end{quotation}

\bigskip

Now, by arguing by induction on $n$ as in \cite[Proof of Theorem 3.7]{BH} 
and by using Lemma \ref{double}, 
we are reduced to the case where $C=(k,l)$ with $k$, $l \ge 1$ and 
$k+l=n$ and $D=(-n)$. Then this follows from \cite[Example 2.25]{BH}.
\end{proof}

\bigskip

\soussection{Some morphisms of algebras ${\boldsymbol{R\Sigma'(W_n) \to R}}$} 
If $\l \in \Bip(n)$, let $\pi_\l : \Sigma'(W_n) \to \ZM$, 
$x \mapsto  \th_n(x)(\cox_\l)$. Recall that 
$\th_n(x)$ is a $\ZM$-linear combination of permutation characters, so 
$\th_n(x)(w)$ lies in $\ZM$. Moreover, $\pi_\l$ does not depend on 
the choice of $\cox_\l$ in $\CC(\l)$, and is a morphism of $\ZM$-algebras. 
We denote by $\pi_\l^R : R\Sigma'(W_n) \to R$ the morphism 
of algebras $\Id_R \otimes_\ZM \pi_\l$. 
We denote by $\DC_\l^R$ the left $R\Sigma'(W_n)$-module whose 
underlying $R$-module is free of rank one and on which $R\Sigma'(W_n)$ 
acts through $\pi_\l^R$. If $K$ is a field, then $\DC_\l^K$ is a simple 
$K\Sigma'(W_n)$-module. 

If $C$ and $D$ are two signed compositions of $n$, let 
$X_{CD}^\subset = \{d \in X_{CD}~|~\lexp{d^{-1}}{W_C} \subset W_D\}$. Then 
\equat\label{valeur caractere}
\pi_{\lamb(C)}(x_D) = |X_{CD}^\subset|.
\endequat

\begin{proof}
By definition, we have
$$\pi_{\lamb(C)}(x_D)=\Bigl(\Ind_{W_D}^{W_n} 1_D\Bigr)(\cox_C)
=\Bigl(\Res_{W_C}^{W_n} \Ind_{W_D}^{W_n} 1_D\Bigr)(\cox_C).$$
Therefore, by the Mackey formula, 
$$\pi_{\lamb(C)}(x_D)=\sum_{d \in X_{CD}} 
\Bigl(\Ind_{W_{C \cap \lexp{d}{D}}}^{W_C} 1_{C \cap \lexp{d}{D}}\Bigr)(\cox_C).$$
But, by the argument in the proof of \cite[proposition 3.12]{BH}, we get 
that $\cox_C$ lies in a subgroup of $W_C$ conjugate to $W_{C \cap \lexp{d}{D}}$ 
if and only if $C \cap \lexp{d}{D}=C$. This shows the result.
\end{proof}

\bigskip

\soussection{Action of the normalizer} 
If $C$ and $D$ are two signed compositions of $n$, we set 
$X_{CD}^\equiv =\{d \in X_{CD}~|~W_C=\lexp{d}{W_D}\}$. Then
\equat\label{wcd}
X_{CD}^\equiv =\{d \in W_n~|~\D_C=d(\D_D)\}.
\endequat

\begin{proof}
Let $\XC=\{d \in W_n~|~\D_C=d(\D_D)\}$. Then it is clear 
that $\XC \subset X_{CD}^\equiv$. Conversely, if $d \in X_{CD}^\equiv$, 
then $d(\Phi_D)=\Phi_C$. So $d(\D_D)$ is a basis of $\Phi_C$, hence 
there exists $w \in W_C$ such that $d(\D_D)=w(\D_C)$. So $d^{-1} w \in X_C$  
and $d^{-1} \in X_C$. So $w=1$, as desired.
\end{proof}

\remark{conjugaison compositions}
If $D$ and $D'$ are two signed compositions of $n$ and if 
$d \in X_{DD'}$ is such that $\lexp{d}{W_{D'}}=W_D$, then 
$d(\D_{D'})=\D_D$ by \ref{wcd} and 
$$X_{D'}=X_D d\qquad \text{and} \qquad x_{D'}=x_D d.$$
Moreover, for every $C \subset D'$, we have 
$d x_C^{D'} d^{-1} = x_{\lexp{d}{C}}^D$ (here, note that 
$d \in X_{DC}$, and we denote for simplification 
$\lexp{d}{C} \cap D$ by $\lexp{d}{C}$ because $\lexp{d}{W_C} \cap W_D = \lexp{d}{W_C}$). 
So conjugacy by $d$ 
induces a morphism of algebras $d_* : R\Sigma'(W_{D'}) \to R\Sigma'(W_D)$.\finl

\bigskip

If $C \equiv D$, then $X_{CD}^\equiv=X_{CD}^\subset$. 
If $D \in \Comp(n)$, we define $\WC(D)=X_{DD}^\subset$. 

\begin{lem}\label{divisible}
Let $D$ be a signed compositions of $n$. Then:
\begin{itemize}
\itemth{a} $\WC(D)=\{w \in W_n~|~w(\D_D)=\D_D\}$.

\itemth{b} $\WC(D)$ is a subgroup of $N_{W_n}(W_D)$.

\itemth{c} The natural map $\WC(D) \to N_{W_n}(W_D)/W_D$ is an isomorphism of groups. 

\itemth{d} $N_{W_n}(W_D) = \WC(D) \ltimes W_D$.

\itemth{e} If $C \in \Comp(n)$, then $|\WC(D)|$ divides $|X_{CD}^\subset|$.
\end{itemize}
\end{lem}

\begin{proof}
(a), (b), (c) and (d) follow immediately from \ref{wcd}. 
Let us now prove (e). First, 
by \ref{valeur caractere}, $|X_{CD}^\subset|$ is equal to the 
number of fixed points of $\cox_C$ in its action on $W_n/W_D$ 
by left multiplication. But $\WC(D)$ acts on 
$W_n/W_D$ by right translation and this action commutes with the left action 
of $\cox_C$. Therefore, $\WC(D)$ permutes the fixed points of $\cox_C$. 
Since $\WC(D)$ acts freely on $W_n/W_D$, (e) follows.
\end{proof}

\bigskip

\section{On the ideals of $R\Sigma'(W_n)$\label{section ideal}}

\medskip

This section is inspired by \cite[\S 3.A]{BP}. 
We shall define some families of left, right and two-sided ideals of $R\Sigma'(W_n)$ 
related to the order $\infspe$ and the preorder $\subset_\lamb$ defined 
in the previous section. We need the following definition: 
if $x \in R\Sigma'(W_n)$, the {\it support of $x$} (denoted by $\Supp(x)$) 
is the subset of $\Comp(n)$ defined by 
$$\Supp(x)=\{C \in \Comp(n)~|~\xi_C(x) \neq 0\}.$$

\bigskip

\soussection{Some left ideals}
A subset $\FC$ of $\Comp(n)$ is called {\it left-saturated} if, for every $D \in \FC$ 
and every $C \in \Comp(n)$ such that $C \infspe D$, we have $C \in \FC$. 
By Proposition C (a), if $\FC$ is left-saturated, then 
$R\Sigma_\FC'(W_n)$ is a left ideal of $R\Sigma'(W_n)$. 

If $x \in R\Sigma'(W_n)$, we set
$$\Sat_l(x)=\{C \in \Comp(n)~|~\exists~D \in \Supp(x),~C \infspe D\}.$$
Then $\Sat_l(x)$ is the minimal left-saturated subset of $\Comp(n)$ 
containing the support of $x$. By the previous remark, 
\equat
R\Sigma'(W_n) x \subset R\Sigma_{\Sat_l(x)}'(W_n).
\endequat

\bigskip

\example{inferieur sature} 
If $D \in \Comp(n)$, then $\FC_{\infspe D}$ and $\FC_{\inferieur D}$ are 
left-saturated. In fact, $\FC_{\infspe D}=\Sat_l(x_D)$. 
Consequently, $R\Sigma_{\infspe D}'(W_n)$ and 
$R\Sigma_{\inferieur D}'(W_n)$ are left ideals of $R\Sigma'(W_n)$. Note that 
$R\Sigma_{\infspe D}'(W_n)/R\Sigma_{\inferieur D}'(W_n)$ is a left 
$R\Sigma'(W_n)$-module which is free of rank $1$ over $R$ (it is generated by 
the image of $x_D$). The action of $R\Sigma'(W_n)$ on this module is 
described in the next proposition.\finl

\begin{prop}\label{tau bis}
Let $D$ be a signed composition of $n$ and let $x \in \Sigma'(W_n)$. Then 
$$xx_D-\pi_{\lamb(D)}^R(x)x_D \in R\Sigma_{\inferieur D}'(W_n).$$
In other words, 
$R\Sigma_{\infspe D}'(W_n)/R\Sigma_{\inferieur D}'(W_n) \simeq \DC_{\lamb(D)}^R$. 
\end{prop}

\begin{proof}
By Proposition C, we only need to show that $\xi_D(xx_D) = \pi_{\lamb(D)}(x)$ 
for every $x \in R\Sigma'(W_n)$. Let $C \in \Comp(n)$. 
By Proposition C, we have 
$$\xi_D(x_C x_D)=|\{d \in X_{CD}~|~W_D \subset \lexp{d^{-1}}{W_C}\}|=|X_{CD}^\subset|.$$
So the result follows from \ref{valeur caractere}. 
\end{proof}

The next result follows immediately from Proposition \ref{tau bis}.

\begin{coro}\label{caractere sature}
Let $\FC$ be a left-saturated subset of $\Comp(n)$ and let 
$\chi_\FC$ denote the character of the left $K\Sigma'(W_n)$-module 
$K\Sigma_\FC'(W_n)$. then 
$$\chi_\FC=\sum_{C \in \FC} \pi_{\lamb(C)}^K.$$
\end{coro}

If $a \in R\Sigma'(W_n)$, we denote by $f_a(T) \in R[T]$ its minimal polynomial. 
Let $\g_a : R\Sigma'(W_n) \to R\Sigma'(W_n)$, 
$x \mapsto ax$ be the left multiplication by $a$ and let $\G_a$ be the matrix of 
$\g_a$ in the basis $(x_C)_{C \in \Comp(n)}$. Then $f_a$ is the minimal 
polynomial of $\g_a$ (or of the matrix $\G_a$). 
By \ref{tau bis}, $\G_a$ is triangular 
(with respect to the order $\infspe$ on $\Comp(n)$) and its characteristic 
polynomial is 
\equat\label{characteristic polynomial}
\prod_{C \in \Comp(n)} (T-\pi_{\lamb(C)}^R(a)).
\endequat
In particular:

\begin{coro}\label{split}
$f_a$ is split over $R$.
\end{coro}

\bigskip

\soussection{Some right ideals}
A subset $\FC$ of $\Comp(n)$ is called {\it right-saturated} if, 
for every $D \in \FC$ and every $C \in \Comp(n)$ such that 
$C \subset_\lamb D$, we have $C \in \FC$. 
If $\FC$ is right-saturated, then $R\Sigma_\FC'(W)$ is a right 
ideal of $R\Sigma'(W_n)$ (see Proposition D (a)). 
It must be noticed that, by opposition with the case of the classical 
Solomon algebra \cite[\S 3.B]{BP}, $R\Sigma_\FC'(W)$ is not necessarily a two-sided 
ideal of $R\Sigma'(W_n)$ (see Example \ref{droite bilatere} below). 

\bigskip

\example{inferieur droite sature} 
If $D \in \Comp(n)$, then $\FC_{\subset_{\lamb} D}$ and $\FC_{\varsubsetneq_\lamb D}$ are 
left-saturated. In fact, $\FC_{\subset_\lamb D}=\Sat_r(x_D)$. 
Consequently, $R\Sigma_{\subset_\lamb D}'(W_n)$ and 
$R\Sigma_{\varsubsetneq_\lamb D}'(W_n)$ are right ideals of $R\Sigma'(W_n)$.\finl

\bigskip

\example{droite bilatere}
Assume here that $n=2$. Then 
$R\Sigma_{\subset_\lamb (-2)}(W_2)=R x_{(-2)}\oplus R x_{(-1,-1)}$ is not a 
two-sided ideal of $\Sigma_R'(W_n)$ because 
$$x_{(1,1)} x_{(-2)}=x_{(-1,-1)} + x_{(-1,1)} - x_{(1,-1)}.~\SS{\square}$$

\bigskip

If $x \in R\Sigma'(W_n)$, we set
$$\Sat_r(x)=\{C \in \Comp(n)~|~\exists~D \in \Supp(a),~C \subset_\lamb D\}.$$
Then $\Sat_r(x)$ is the minimal right-saturated subset of $\Comp(n)$ 
containing the support of $x$. By the previous remark, 
\equat\label{inclusion droite}
xR\Sigma'(W_n) \subset R\Sigma_{\Sat_r(x)}'(W_n).
\endequat
We shall construct in Section \ref{section positif} a class 
of elements $x$ for which equality holds in \ref{inclusion droite}. 

\bigskip

\soussection{Some two-sided ideals\label{section saturee}}
A subset $\FC$ of $\Comp(n)$ is called {\it saturated} if it is 
left-saturated and right-saturated. If $\FC$ is saturated, then 
$R\Sigma_\FC'(W_n)$ is a two-sided ideal of $R\Sigma'(W_n)$. 

\bigskip

\example{exemple bilateres} 
If $k \ge 0$, we set $\FC_k(n)=\{C \in \Comp(n)~|~\length(C) \ge k+1\}$ 
and $\FC_k^-(n)=\{C \in \Comp(n)~|~\length(C) + \length^-(C) \ge k+1\}$. 
Then, by the Remarks \ref{proprietes inclusion} and \ref{comparaison}, 
$\FC_k(n)$ and $\FC_k^-(n)$ are saturated subsets of $\Comp(n)$.\finl

\bigskip 

\section{Positivity properties\label{section positif}}

\medskip

In this section, and only in this section, we assume that $K$ is an ordered 
(for instance $K=\QM$ or $K=\RM$). Recall that this implies 
that $K$ has characteristic $0$. 
We shall now construct a class of elements of $K\Sigma'(W_n)$ for which 
equality holds in \ref{inclusion droite}. 
We denote by $K\Sigma'(W_n)^+$ the set of elements $a \in K\Sigma'(W_n)$ such that 
$\xi_C(a) \ge 0$ for every $C \in \Comp(n)$. Note that $x_C \in K\Sigma'(W_n)^+$ 
for any $C \in \Comp(n)$. If $a$ and $b$ are two elements of $K\Sigma'(W_n)^+$, then 
\equat
a+b \in K\Sigma'(W_n)^+.
\endequat
However, contrarily to the case of Solomon algebras \cite[3.2]{BP}, it might happen 
that $ab \not\in K\Sigma'(W_n)^+$ (see example \ref{droite bilatere}). 
However, the analogue of \cite[First statement of Proposition 3.6]{BP} holds:

\begin{prop}\label{sature droite}
Assume that $K$ is an ordered field. 
Let $a \in K\Sigma'(W_n)^+$. Then
$$aK\Sigma'(W_n)=K\Sigma_{\Sat_r(a)}'(W_n).$$
\end{prop}

\begin{proof}
Let $\FC=\Sat_r(a)$. By \ref{inclusion droite}, we have 
$aK\Sigma'(W_n) \subset K\Sigma_{\FC}'(W_n)$. We shall show by induction on 
$C \in \FC$ (induction with respect to the order $\infspe$) that 
$x_C \in a K\Sigma'(W_n)$. For this, we may, and we will, assume that 
$a \not= 0$. 

First, if $C=(-1,-1,\dots,-1)$, then $x_C=\sum_{w \in W_n} w$ so 
$$ax_C=\Bigl(\sum_{D \in \Comp(n)} \xi_D(a) |X_D| \Bigr) x_C.$$
Since $a \neq 0$ and $a \in K\Sigma'(W_n)^+$, 
we have by definition $\sum_{D \in \Comp(n)} \xi_D(a) |X_D| > 0$. So 
$x_{(-1,-1,\dots,-1)} \in a K\Sigma'(W_n)$. 

Now, let $C \in \FC$ and assume that, if $C' \in \FC$ is such that $C' \inferieur C$, 
then $x_{C'} \in a K\Sigma'(W_n)$. Then, by Propositions C and \ref{tau bis}, 
we have 
$$a x_C -\pi_{\lamb(C)}(a) x_C \in K\Sigma_{\inferieur C}'(W_n).$$
But, by the induction hypothesis, we have that 
$K\Sigma_{\subset_\lamb C}'(W_n) \subset a K\Sigma'(W_n)$. So it remains 
to show that $\pi_{\lamb(C)}(a) \neq 0$. But,
$$\pi_{\lamb(C)}(a) = \sum_{D \in \Supp(a)} \xi_D(a) \pi_{\lamb(C)}(x_D).$$
Since $\xi_D(a) > 0$ and $\pi_{\lamb(C)}(x_D) \ge 0$ for every $D \in \Supp(a)$, 
it remains to show that there exists $D \in \Supp(a)$ such that 
$\pi_{\lamb(C)}(x_D) > 0$. But, by the definition of $\FC$, there exists 
$D \in \Supp(a)$ such that $W_C$ is contained in a conjugate of $W_D$. 
So, for this particular $D$, we have that $\cox_C$ is contained in a conjugate 
of $W_D$. So $\pi_{\lamb(C)}(x_D)=\Ind_{W_D}^{W_n}(\cox_C) \ge 1$ 
and the proof of the proposition is complete.
\end{proof}

The next four corollaries must be compared with 
\cite[Corollaries 4.7, 3.8 and 3.12 and Proposition 3.10]{BP}. 

\begin{coro}\label{inversible}
Assume that $K$ is an ordered field. Let $a \in K\Sigma'(W_n)^+$. 
Then $a$ is invertible in $K\Sigma'(W_n)$ 
if and only if $\xi_n(a) > 0$.
\end{coro}

\begin{coro}\label{somme ideaux}
Assume that $K$ is an ordered field. 
Let $a_1$,\dots, $a_r \in K\Sigma'(W_n)^+$. Then 
$a_1+\dots + a_r \in K\Sigma'(W_n)^+$ and 
$$a_1 K\Sigma'(W_n) + \dots + a_r K\Sigma'(W_n)= 
(a_1+\dots+a_r) K\Sigma'(W_n).$$
\end{coro}

The proof of the next corollary follows an argument of Atkinson \cite{atkinson}.

\begin{coro}\label{minimal positif}
Assume that $K$ is an ordered field. 
Let $a \in K\Sigma'(W_n)^+$ and let $r$ be a non-zero natural number. Then:
\begin{itemize}
\itemth{a} $\val f_a \le 1$.

\itemth{b} $a^r K\Sigma'(W_n)=aK\Sigma'(W_n)$.

\itemth{c} $K\Sigma'(W_n)a^r=K\Sigma'(W_n)a$.
\end{itemize}
\end{coro}

\begin{proof}
Recall that $f_a$ is the minimal polynomial of $a$. We first 
prove (b). It is sufficient to show the result for $r=2$. 
Let $m : aK\Sigma'(W_n) \to aK\Sigma'(W_n)$, $x \mapsto ax$. 
Then, by Proposition \ref{sature droite}, we have 
$aK\Sigma'(W_n)=K\Sigma_{\Sat_r(a)}'(W_n)$. But, by the proof of 
Proposition \ref{sature droite}, we have $\pi_{\lamb(C)}(a) > 0$ for every 
$C \in \Sat_r(a)$. Therefore, by Proposition \ref{tau bis}, 
the matrix of $m$ in the basis $(x_C)_{C \in \Sat_r(a)}$ 
is triangular (with respect to the order $\infspe$) and has 
positive diagonal coefficients. So it is invertible. This shows (b). 

\medskip

(a) By (b), the minimal polynomial $f \in K[T]$ of $m$ has a non-zero 
constant term. But, $f(a)a=f(m)(a) = 0$. Therefore, $f_a$ divides 
$T f(T)$. This shows (a).

\medskip

(c) Now, by (a), we have that $a \in K[a] a^2$. So $a \in K\Sigma'(W_n) a^2$, 
as desired.
\end{proof}

Recall that $\g_a$ denote the left multiplication 
$K\Sigma'(W_n) \to K\Sigma'(W_n)$, $x \mapsto ax$. 
Let $\d_a : K\Sigma'(W_n) \to K\Sigma'(W_n)$, $x \mapsto xa$ denote the right 
multiplication by $a$. 

\begin{coro}\label{image noyau}
Assume that $K$ is an ordered field. 
Let $a \in K\Sigma'(W_n)^+$. Then:
\begin{itemize}
\itemth{a} $\Ker \g_a \oplus \im \g_a = K\Sigma'(W_n)$.

\itemth{b} $\Ker \d_a \oplus \im \d_a = K\Sigma'(W_n)$.
\end{itemize}
\end{coro}

\begin{proof}
(a) For dimension reasons, it is sufficient to prove that 
$\Ker \g_a \cap \im \g_a=0$. Let $x \in \Ker \g_a \cap \im \g_a$. 
Then $ax=0$ and there exists $y \in K\Sigma'(W_n)$ such that $x=ay$. 
So $a^2 y = 0$. Therefore, $a^ry=0$ for every $r \ge 2$. 
But $a \in \sum_{r \ge 2} K a^r$ by Corollary \ref{minimal positif} (a). 
So $ay=0$. In other words, $x=0$, as desired. The proof of (b) is similar.
\end{proof}

\remark{pas semisimple}
By opposition with the case of Solomon descent algebras, it may happen 
that $f_{x_C}$ is not square-free (compare with \cite[Proposition 3.10]{BP}). For 
instance, if $n=4$ and if $K=\QM$, we have 
$$f_{x_{(-3,1)}}(T)=T(T-2)(T-4)(T-8)^2(T-32).$$
This computation has been done using {\tt CHEVIE} \cite{chevie}.\finl

\bigskip

We close this section by a result on the centralizers of positive elements 
(compare with \cite[Corollary 3.12]{BP}: the proof presented 
here is really different):

\begin{prop}\label{prop centralisateur}
Assume that $K$ is an ordered field. 
Let $a \in K\Sigma'(W_n)^+$ and $r$ be a non-zero natural number. Then 
$Z_{K\Sigma'(W_n)}(a)=Z_{K\Sigma'(W_n)}(a^r)$. 
\end{prop}

\begin{proof}
Let $A=\End_K K\Sigma'(W_n)$. Let $\g : K\Sigma'(W_n) \to A$, 
$x \mapsto \g_x$. It is an injective homomorphism of algebras. Therefore, 
$Z_{K\Sigma'(W_n)}(a)=\g^{-1}(Z_A(\g_a))$. So, in order to prove 
the proposition, we only need to prove that $Z_A(\g_a)=Z_A(\g_a^r)$. 

Let $A'=\End_K(\Ker \g_a)$ and $A''=\End_K(\im \g_a)$. Then, by 
Corollary \ref{image noyau},  
$A' \oplus A''$ is a sub-$K$-algebra of $A$ and $Z_A(\g_a)$ 
is contained in $A' \oplus A''$. Let $\g''$ denote the restriction 
of $\g_a$ to $\im \g_a$. Then $Z_A(\g_a)=A' \oplus Z_{A''}(\g'')$. 
Since $\Ker \g_a^r = \Ker \g_a$ and $\im \g_a^r=\im \g_a$, 
we only need to prove that $Z_{A''}(\g'')=Z_{A''}(\g^{\prime\prime r})$. 
But, by Proposition \ref{sature droite} and its proof, 
$(x_C)_{C \in \Sat_r(a)}$ is a basis of $\im \g_a$ and 
the matrix of $\g''$ in this basis is triangular 
(with respect to the order $\infspe$) with positive coefficients 
on the diagonal. So the proposition follows from Lemma 
\ref{centralisateur triangulaire} below.
\end{proof}

\begin{quotation}
{\small 
\begin{lem}\label{centralisateur triangulaire}
Let $m$ be a non-zero natural number. Let $M=(m_{ij}) \in \Mat_m(K)$ 
be an upper triangular $m \times m$ matrix such that $m_{ii} > 0$ 
for every $i \in \{1,2,\dots,m\}$. Then 
$$Z_{\Mat_m(K)}(M)=Z_{\Mat_m(K)}(M^r)$$
for every $r \ge 1$.
\end{lem}

\begin{proof}[Proof of Lemma \ref{centralisateur triangulaire}]
Let $E=\Mat_m(K)$. Since $M$ is invertible, we can write $M=SU=US$ where 
$S$ (resp. $U$) is a diagonalizable (resp. unipotent) matrix. 
Then $Z_E(M)=Z_E(S) \cap Z_E(U)$. So it is sufficient to show that 
$Z_E(S)=Z_E(S^r)$ and $Z_E(U)=Z_E(U^r)$. 

Since $S$ is diagonalizable, we may assume that it is diagonal. 
Now the fact that $Z_E(S)=Z_E(S^r)$ follows from the fact that, 
if $x$, $y \in K$ are such that $x> 0$, $y> 0$ and $x^r=y^r$, then $x=y$ 
(because $K$ is an ordered field).

Let us now show that $Z_E(U)=Z_E(U^r)$. Since $K$ is an ordered field, 
its characteristic is zero. Let $N$ be a nilpotent matrix such that 
$U=e^N$ (exponential). Then $Z_E(U)=Z_E(N)$ and, since $U^r=e^{rN}$, 
we have $Z_E(U^r)=Z_E(rN)=Z_E(N)$ (because $r \neq 0$ in $K$).
\end{proof}}
\end{quotation}

\bigskip

\section{Action of the longest element\label{action plus long}}

\medskip

If $C \in \Comp(n)$, we denote by $w_C$ the longest element of $W_C$. 
If $C \in \Comp^+(n)$, we denote by $\s_C$ the longest element of $\SG_C$ 
(in other words, $\s_C=w_{-C}$). 
In particular, $w_n$ is the longest element of $W_n$. Recall that 
$w_n \in R\Sigma'(W_n)$ (in fact $w_n \in R\Sigma(W_n)$) 
and that $\th_n^R(w_n)=\e_n$, the sign character of $W_n$ (see for instance 
\cite{solomon}). Moreover, $w_n$ is central in $W_n$, so it is central 
in $R\Sigma'(W_n)$. 

\medskip

First, note that 
\equat\label{signe coxeter}
\e_n(\cox C)=(-1)^{n-\length^-(C)}.
\endequat
Recall that the function $\length^- : \Comp(n) \to \NM$ has been 
defined in \S\ref{section signed}. In particular, 
\equat\label{tau signe}
\pi_\l(w_n)=(-1)^{n - \length^-(\l)}
\endequat
for all $\l \in \Bip(n)$. 

\bigskip

{\it From now on, and until the end of this section, we assume 
that $2$ is invertible in $R$.} 
Write
$$e_n^+=\frac{1}{2}(1+w_n)\qquad\text{and}\qquad e_n^-=\frac{1}{2}(1-w_n).$$
Since $w_n$ is central in $W_n$, $e_n^+$ and $e_n^-$ are central idempotents 
of $K\Sigma'(W_n)$. Moreover, they are orthogonal and $e_n^+ + e_n^-=1$. 

We shall now describe a basis of $R\Sigma'(W_n)e_n^+$ and 
$R\Sigma'(W_n)e_n^-$. 
This is build on the same model as for the classical Solomon algebra 
\cite[\S 2.B]{BP}. First, note that 
\equat\label{trivial wn}
w_n X_C = X_C w_C
\endequat
for every signed composition $C$ of $n$. This is proved as follows. 
An element $w \in W_n$ belongs to $X_C$ if and only if $w(\a) > 0$ 
for every $\a \in \Phi_C^+$. 
In other words, $w \in w_nX_C$ if and only if $w(\a) < 0$ for every 
$\a \in \Phi_C^+$. A similar argument applies to show that $w \in X_C w_C$ 
if and only if $w(\a) < 0$ for every $\a \in \Phi_C^+$. This shows \ref{trivial wn}. 
Consequently, 
\equat\label{trivial wn bis}
w_n x_C = x_C w_C.
\endequat

Since $2$ is invertible in $R$, we can define 
$$x_C'=\sum_{\substack{D \in \Comp(n) \\ S_D \subset S_C}} 
\Bigl(-\frac{1}{2}\Bigr)^{|S_C|-|S_D|} x_D.$$
These elements have the following properties:

\begin{prop}\label{base xprime}
If $2$ is invertible in $R$, then:
\begin{itemize}
\itemth{a} $(x_C')_{C \in \Comp(n)}$ is an $R$-basis of $R\Sigma'(W_n)$.

\itemth{b} $\Ker \th_n^R=\sum_{C \equiv D} R(x_C'-x_D')$.

\itemth{c} $w_n x_C'=(-1)^{|S_C|} x_C'=(-1)^{n-\length^-(C)} x_C'$.
\end{itemize}
\end{prop}

\begin{proof}
(a) is trivial. Note for information that 
$$x_C=\sum_{\substack{D \in \Comp(n) \\ S_D \subset S_C}} 
\Bigl(\frac{1}{2}\Bigr)^{|S_C|-|S_D|} x_D'.$$

(b) follows from Theorem B (c) and Remark \ref{conjugaison compositions}. 

Let us now prove (c). First, we set 
$$\xti_C=\sum_{\substack{D \in \Comp(n) \\ S_D \subset S_C}} 
\Bigl(-\frac{1}{2}\Bigr)^{|S_C|-|S_D|} x_D^C.$$
Now, by \ref{trivial wn bis}, we get 
$w_n x_C = x_C w_C \xti_C$. But now, $\xti_C$ is an element of the classical 
Solomon descent algebra $R\Sigma(W_C)$ to which the result of \cite[2.12]{BP} 
can be applied: we get $w_C \xti_C =(-1)^{|S_C|} \xti_C$. This shows the 
first equality. The second equality is easy from the definition of $S_C$.
\end{proof}

\begin{coro}\label{dim paire}
If $2$ is invertible in $R$, then 
$\dim_R R\Sigma'(W_n) e_n^+=\dim_R R\Sigma'(W_n) e_n^- = 3^{n-1}$.
\end{coro}

\begin{proof}
Let 
$$\CC^+=\{C \in \Comp(n)~|~\length^-(C) \equiv n \mod 2\}$$
$$\CC^-=\{C \in \Comp(n)~|~\length^-(C) \equiv n+1 \mod 2\}.\leqno{\text{and}}$$ 
Then, by Proposition \ref{base xprime} (a) and (c), we have 
$\dim_R R\Sigma'(W_n)e_n^? = |\CC^?|$ for $? \in \{+,-\}$. 
It is now sufficient to show that $|\CC^+|=|\CC^-|=3^{n-1}$. Since 
the number of compositions of $n$ of length $k$ is equal to 
$\begin{pmatrix} n-1 \\ k-1\end{pmatrix}$, 
we have 
$$|\CC^+|=\sum_{k=1}^n \begin{pmatrix} n-1 \\ k-1 \end{pmatrix} 2^{k-1}=3^{n-1},$$
and similarly for $|\CC^-|$.
\end{proof}

\example{x un un prime}
We consider in this example the elements of the form 
$x_{\nu_1,\nu_2,\dots,\nu_n}'$ of $R\Sigma'(W_n)$ where $\nu_i \in \{1,-1\}$. 
We shall show that they are quasi-idempotents. We first need some notation. 
Let $I_n^+=\{1,2,\dots,n\}$. 
If $\s \in W_n$, let $\Ib(\s)=\{i \in I_n^+~|~\s(i) > 0\}$. 
If $I \subset I_n^+$, we denote by $\Cb(I)$ the signed composition 
$(\nu_1,\dots,\nu_n)$, where $\nu_i=1$ (respectively 
$\nu_i=-1$) if $i \in I$ (respectively $i \not\in I$). 
We also define $\g_I : W_n \to R^\times$, $\s \mapsto (-1)^{|I|-|I \cap \Ib(\s)|}$. 
For instance, $\g_\vide=1_n$ and $\g_{I_n^+}$ is a linear character 
of $W_n$ (it will be denoted by $\g_n$ for simplification). Note also 
that the restriction of $\g_I$ to $\TG_n$ is always a linear character 
(the group $\TG_n$ has been defined in \S\ref{section type B}: 
it is the group generated by $t_1$,\dots, $t_n$). Moreover, 
if $\s \in \SG_n$ and $\t \in W_n$, we have $\Ib(\s\t)=\Ib(\t)$, 
so $\g_I(\s \t)=\g_I(\t)$. Finally, we denote by $\SG_n(I)$ the 
stabilizer of $I$ in $\SG_n$. 

With this notation, we have
\equat\label{x un un}
x_{\Cb(I)}' = \frac{1}{2^{|I|}} \sum_{\s \in W_n} \g_I(\s) \s.
\endequat
\begin{quotation}
\begin{proof}[Proof of \ref{x un un}]
First, note that 
$$x_{\Cb(I)}'= \sum_{J \subset I} 
\Bigl(-\frac{1}{2}\Bigr)^{|I|-|J|} x_{\Cb(J)}.$$

Now, let $\s \in W_n$ and $J \subset I_n^+$. 
Then $\s \in X_{\Cb(J)}$ if and only $J \subset \Ib(\s)$. 
Therefore, by the previous equality, the coefficient of $\s$ in $x_{\Cb(I)}'$ is equal to 
$$\sum_{J \subset I \cap \Ib(\s)} \Bigl(-\frac{1}{2}\Bigr)^{|I|-|J|}=
\frac{(-1)^{|I|-|I \cap \Ib(\s)|}}{2^{|I|}},$$
as desired. 
\end{proof}
\end{quotation}

If $I$ and $J$ are two subsets of $I_n^+$, then 
\equat\label{x un un quasi}
x_{\Cb(I)}'x_{\Cb(J)}' = \begin{cases}
2^{n-|I|} |\SG_n(I)| x_{\Cb(J)}' & \text{if $|I|=|J|$},\\
0 & \text{otherwise.}
\end{cases}
\endequat

\begin{quotation}
\begin{proof}[Proof of \ref{x un un quasi}]
Let $e_I=\sum_{\t \in \TG_n} \g_I(\t) \t$ and let 
$e=\sum_{\s \in \SG_n} \s$. Then, by \ref{x un un}, we have 
$x_{\Cb(I)}'=ee_I/2^{|I|}$. 
Moreover,  
$$e_I e_J =\begin{cases} 2^n e_I & \text{if $I=J$},\\ 0 & \text{otherwise.}\end{cases}$$
Therefore,
\eqna
x_{\Cb(I)}'x_{\Cb(J)}'&=&\DS{\frac{1}{2^{|I|+|J|}} e e_I \Bigl(\sum_{\s \in \SG_n} \s\Bigr)
e_J }\\
&=& \DS{\frac{1}{2^{|I|+|J|}} \sum_{\s \in \SG_n} e \s^{-1} e_I \s e_J }\\
&=& \DS{\frac{1}{2^{|I|+|J|}} \sum_{\s \in \SG_n} e e_{\s^{-1}(I)} e_J }.\\
\endeqna
If $|I| \neq |J|$, then $\s^{-1}(I) \neq J$ for every $\s \in \SG_n$. 
If $|I|=|J|$, then the number of elements $\s \in \SG_n$ such that 
$\s^{-1}(I)=J$ is equal to $|\SG_n(I)|=|I|!(n-|I|)!$. This shows 
the result in this last case. 
\end{proof}
\end{quotation}

In particular, if $2 |\SG_n(I)|$ is invertible in $R$, then 
$x_{\Cb(I)}'/(2^{n-|I|} |\SG_n(I)|)$ is an idempotent of $R\Sigma'(W_n)$. 
We shall show later that it is a primitive idempotent of 
$R\Sigma'(W_n)$ (see \ref{conjugue x un un}). A description of 
the module $\QM W_n x_{\Cb(I)}'$ will be given in Example 
\ref{exemple question 6}.

\medskip

Since $\g_{I_n^+}=\g_n$ is a linear character of $W_n$, we deduce 
immediately the following two properties of $x_{1,1,\dots,1}'$:
\equat\label{central x un}
\text{\it $x_{1,1,\dots,1}'$ is central in $\QM W_n$, hence is central in $\QM\Sigma'(W_n)$;}
\endequat
\equat\label{x un idem}
(x_{1,1,\dots,1}')^2 = n!~ x_{1,1,\dots,1}'.
\endequat
In particular, if $p$ does not divide $|W_n|$, then $x_{1,1,\dots,1}/n!$ is a primitive 
central idempotent of $\QM W_n$, hence is a primitive central idempotent 
of $\QM\Sigma'(W_n)$.\finl

\bigskip

\section{Restriction morphisms between Mantaci-Reutenauer algebras}

\medskip

F. Bergeron, N. Bergeron, R.B. Howlett and D.E. Taylor \cite{bbht} 
have constructed so-called {\it restriction morphisms} 
between the Solomon algebra of a finite Coxeter group and the Solomon 
algebras of its standard parabolic subgroups. We shall construct 
here a restriction morphism $R\Sigma'(W_n) \to R\Sigma'(W_D)$ whenever 
$D$ is a semi-positive signed compositions of $n$. It might be possible 
that such a morphism exists for every signed compositions, but we are not 
able to prove it (or to prove that there is no analogue in general). 
Most of the results of this section have an analogue in the context 
of Solomon's descent algebras \cite[\S 4]{BP}.

\bigskip

\soussection{Definition\label{def restriction}} 
We fix in this section a semi-positive signed composition $D$ of $n$.
We denote by $\Res_D : R\Sigma'(W_n) \to R\Sigma'(W_D)$ the unique 
$R$-linear map such that 
$$\Res_D x_C=\sum_{d \in X_{CD}} x_{\lexp{d^{-1}}{C} \cap D}^D$$
for every $C \in \Comp(n)$. If $C \subset D$ is a semi-positive 
signed composition, we define $\Res_C^D : R\Sigma'(W_D) \to R\Sigma'(W_C)$ 
similarly. 

\begin{prop}\label{restriction}
Let $D$ be a semi-positive signed composition of $n$. Then:
\begin{itemize}
\itemth{a} If $x \in \Sigma(W_n)$, then $x_D \Res_D(x) = xx_D$.

\itemth{b} $\Res_D$ is a morphism of algebras.

\itemth{c} If $C \subset D$ is also semi-positive, then 
$\Res_C^D \circ \Res_D = \Res_C$.

\itemth{d} The diagram
$$\diagram
R\Sigma'(W_n)\ddto_{\DS{\Res_D}} \rrto^{\DS{\th_n^R}} && 
R\Irr W_n \ddto^{\DS{\Res_{W_D}^{W_n}}}\\
&&\\
R\Sigma'(W_D) \rrto^{\DS{\th_D^R}} && R\Irr W_D
\enddiagram$$
is commutative.

\itemth{e} If $D'$ is another signed composition of $n$ and if $d \in X_{DD'}$ 
is such that $\lexp{d}{S_{D'}}=S_D$, then $D'$ is semi-positive and 
$d_* \circ \Res_{D'} = \Res_D$.
\end{itemize}
\end{prop}

\begin{proof}
(a) follows from Proposition C (b). (b) and (c) follow from (a) and from the fact 
that the map $\mu_D : R\Sigma'(W_D) \to R\Sigma'(W_n)$, $x \mapsto x_D x$ 
is injective. (d) follows from the Mackey formula. (e) follows easily 
from Remark \ref{conjugaison compositions} and from (a).
\end{proof}

By Remark \ref{conjugaison compositions}, 
the group $\WC(D)$ acts on $R\Sigma'(W_D)$. 
Moreover, by Proposition \ref{restriction} (e), 
we have
\equat\label{inclusion image}
\im \Res_D \subset R\Sigma'(W_D)^{\WC(D)}.
\endequat

Let $\Comp(D)=\{C \in \Comp(n)~|~C \subset D\}$. Write $D=(d_1,\dots,d_r)$ and 
$\Bip(D)=\Bip(d_1) \times \dots \times \Bip(d_r)$: here, if $d < 0$, $\Bip(d)$ denotes 
the set of partitions of $-d$. If $C \in \Comp(D)$, we denote by 
$\lamb_D(C)$ the element of $\Bip(D)$ defined in the natural way 
component by component. 
Then $\lamb_D(C) =\lamb_D(C')$ if and only if $C \equiv_D C'$. Therefore, 
the canonical injection $\Comp(D) \injto \Comp(n)$ induces a unique 
map $\t_D : \Bip(D) \to \Bip(n)$ such that $\t_D(\lamb_D(C))=\lamb(C)$.

\begin{coro}\label{tau compose}
Let $D$ be a semi-positive signed composition of $n$ and let $\l \in \L_D$. 
Then $\pi_{\t_D(\l)}^R = \pi_\l^R \circ \Res_D$.
\end{coro}

\begin{proof}
Let $(\xi_C^D)_{C \in \Comp(D)}$ denote the basis of $\Hom_R(R\Sigma'(W_D),R)$ 
dual to $(x_C^D)_{C \in \Comp(D)}$. 
Let $C$ be a signed composition of $n$ which is contained in $D$ and 
let $\l=\lamb_D(C) \in \Bip(D)$. Then, if $x \in R\Sigma'(W_D)$, 
we have by Proposition \ref{tau bis}, 
$$\pi_\l^R(x)=\xi_C^D(xx_C^D).$$
Therefore, if $x \in R\Sigma'(W_n)$, we have 
\eqna
\pi_\l^R(\Res_D x)&=&\xi_C^D((\Res_D x)x_C^D)\\
                 &=&\xi_C(x_D(\Res_D x)x_C^D) \\
                 &=&\xi_C(xx_Dx_C^D) \\
                 &=&\xi_C(xx_C) \\
                 &=&\pi_{\lamb(C)}(x).
\endeqna
Now, the result follows from the fact that $\lamb(C)=\t_D(\l)$ 
by definition.
\end{proof}

\medskip

We conclude this subsection by a result on the kernel of $\Res_D$.

\begin{prop}\label{noyau}
We have:
\begin{itemize}
\itemth{a} $\Ker(\Res_D) = \{x \in R\Sigma'(W_n)~|~xx_D=0\}$.

\itemth{b} If $K$ is a field of characteristic zero, then 
$K\Sigma'(W_D)=\Ker(\Res_D) \oplus K\Sigma'(W_n)x_D$.
\end{itemize}
\end{prop}

\begin{proof}
(a) follows from Proposition \ref{restriction} (a). To prove (b), we may, 
and we will, assume that $K=\QM$. Then, since $\QM$ is an ordered field, 
(b) follows from (a) and Corollary \ref{image noyau}.
\end{proof}

\bigskip

\soussection{Restriction to the Solomon algebra of ${\boldsymbol{\SG_n}}$} 
If $D$ is a semi-positive signed composition of $n$, then $\SG_{D^+}$ is a 
parabolic subgroup of $\SG_n$. So there is a restriction morphism 
$\Res_{D^+}^\SG : R\Sigma(\SG_n) \to R\Sigma(\SG_{D^+})$ which was constructed 
in \cite{bbht} (see also \cite[Proposition 4.1]{BP} for the proof of the fact that 
it is a morphism of algebras: it works as in the Proposition \ref{restriction} 
above). Then the diagram 
\equat\label{restriction sn}
\diagram
R\Sigma(\SG_n) \xto[0,2]|<\ahook \ddto_{\DS{\Res_{D^+}^\SG}} && 
R\Sigma'(W_n) \ddto^{\DS{\Res_D}} \\
&& \\
R\Sigma(\SG_{D^+}) \xto[0,2]|<\ahook && R\Sigma'(W_D) \\
\enddiagram
\endequat
is commutative. Indeed, if $x \in R\Sigma(\SG_n)$, then, by definition, we have 
$$x_D \Res_D(x)=x x_D=x x_{D^+} x_D^{D^+} = x_{D^+} \Res_{D^+}^\SG(x) x_D^{D^+}.$$
So it remains to show that, if $u \in R\Sigma(\SG_{D^+})$, then 
$u x^{D^+}_D = x_D^{D^+} u$. By direct product, we are reduce 
to prove this whenever $D=D^+$ (in which case it is trivial) or 
whenever $D=D^-$. In this last case, since $D$ is semi-positive, 
we have $D^+=(1)$ and the result follows from the fact 
that the algebra $R\Sigma'(W_1)=RW_1$ is commutative.

\bigskip

\soussection{An example} 
Now, let us consider a particular semi-positive signed composition. 
If $D=(k,-1,-1,\dots,-1)$, where $k \ge 1$, then $\Res_D$ induces 
in fact a morphism of algebras $\Res_k^n : R\Sigma'(W_n) \to R\Sigma'(W_k)$. 
Since $D$ is also parabolic, we have that 
$\Res_k^n(R\Sigma(W_n))\subset R\Sigma(W_k)$, that the induced map 
$R\Sigma(W_n) \to R\Sigma(W_k)$ coincides with the map 
denoted by $\Res_{S_k}^{S_n}$ in \cite[\S 4.B]{BP} and that the diagram
\equat\label{restriction solomon}
\diagram
R\Sigma(W_n) \xto[0,2]|<\ahook \ddto_{\Res_{S_k}^{S_n}} && R\Sigma'(W_n) \ddto^{\Res_k^n} \\
&& \\
R\Sigma(W_k) \xto[0,2]|<\ahook && R\Sigma'(W_k) \\
\enddiagram
\endequat
is commutative. The next result can be compared with \cite[Proposition 4.15]{BP}. 

\begin{prop}\label{surjectif B}
If $R=K$ is a field of characteristic zero, then $\Res_k^n$ is surjective.
\end{prop}

\remark{surjectif Z} 
It is probable that the above proposition remains valid for every 
commutative ring $R$ (i.e. for $R=\ZM$). It has been checked 
for $n \le 5$ using {\tt CHEVIE} \cite{chevie}.\finl

\begin{proof}
By transitivity (see Proposition \ref{restriction} (c)), we only need 
to prove that $\Res_{n-1}^n$ is surjective. We almost reproduce 
the argument in \cite[Proposition 4.15]{BP}. We have
\begin{multline*}
X_{n-1,-1} = \{ s_i s_{i+1} \dots s_{n-1}~|~1 \le i \le n\} \\
\coprod \quad\{s_i s_{i-1} \dots s_1 ts_1s_2 \dots s_{n-1}~|~0 \le i \le n-1\}.
\end{multline*}
Therefore, if $d \in W_n$ and $i \in \{1,2,\dots,n-1\}$ (resp. $i \in \{1,2,\dots,n\}$) 
are such that $d^{-1} \in X_{n-1,-1}$, $\ell(ds_i) > \ell(d)$ (resp. $\ell(dt_i) > \ell(d)$), 
and $d s_i d^{-1} \in S_{n-1}'$ (resp. $d t_i d^{-1} \in S_{n-1}'$), then 
$$ds_i d^{-1} \in \{s_i,s_{i-1}\}\leqno{(*)}$$
(resp. 
$$dt_i d^{-1} \in \{t_i,t_{i-1}\}.\leqno{(**)}$$

We now define a total order $\trianglelefteq$ on $\Comp(n-1)$. 
Let $C$ and $D$ be two signed compositions of $n-1$. We write 
$C \trianglelefteq D$ 
if and only if one of the following two conditions are satisfied:
\begin{quotation}
\begin{itemize}
\itemth{1} $|S_C'| < |S_D'|$;

\itemth{2} $|S_C'|=|S_D'|$ and $S_C'$ is smaller than $S_D'$ for the lexicographic 
order induced by the order $t_1 < t_2 < \dots < t_{n-1} < s_1 < \dots < s_{n-2}$ 
on $S'_{n-1}$. 
\end{itemize}
\end{quotation}
It follows immediately from $(*)$ and $(**)$ that 
$$\Res_{n-1}^n x_{D \sqcup (-1)} \in \a_D x_D + 
\sum_{C \vartriangleleft D} K x_C$$
with $\a_D \in \ZM$, $\a_D > 0$ (for every $D \in \Comp(n-1)$). 
Recall that the operation $\sqcup$ (concatenation) has been 
defined in \S\ref{section signed}. 
The proof of the proposition is complete.
\end{proof}

\bigskip

\section{Simple modules, radical}

\medskip

\begin{quotation}
\noindent{\bf Hypothesis and notation:} {\it 
From now on, and until the end of this paper, 
we assume that $R=K$ is a field. We denote by 
$p$ its characteristic ($p \ge 0$). We denote by $\Bip_{p'}(n)$ the 
set of bipartitions $\l$ of $n$ such that $o(\l)$ is invertible in $K$ 
(recall that $o(\l)$ denotes the order of $\cox_\l$). If $\l \in \Bip(n)$, 
we denote by $\l_{p'}$ the bipartition of $n$ such that the $p'$-part 
of $\cox_{\l}$ is conjugate to $\cox_{\l_{p'}}$ (if $p=0$, then $\l_{p'}=\l$).}
\end{quotation}

\bigskip

\soussection{Simple modules} 
Since $\Ker \th_n^K$ is a nilpotent two-sided ideal of $K\Sigma'(W_n)$, it 
is contained in the kernel of every simple representation of 
$K\Sigma'(W_n)$. Therefore, every simple representation factorizes 
(through $\th_n^K$) to a simple representation of $K\Irr W_n$. 
Since every irreducible character of $W_n$ has value in $\ZM$, 
and since $\th_n^K$ is surjective by Theorem B (c), 
we get (see for instance \cite[Proposition 2.14 and Corollary 2.15]{bonnafe gr}):

\begin{prop}\label{simples}
Let $\l$ and $\mu$ be two bipartitions of $n$. 
\begin{itemize}
\itemth{a} $\DC_\l^K \simeq \DC_\mu^K$ if and only if $\l_{p'}=\mu_{p'}$. 

\itemth{b} $\{\DC_\l^K~|~\l \in \Bip_{p'}(n)\}$ is a set of 
representatives of isomorphy classes of simple left $K\Sigma'(W_n)$-modules.
\end{itemize}
\end{prop}

\begin{coro}
$K\Sigma'(W_n)$ is split.
\end{coro}

\begin{coro}\label{irreductibles}
$\Irr_K K\Sigma'(W_n)=\{\pi_\l^K~|~\l \in \Bip_{p'}(n)\}$. 
\end{coro}

The formula for the irreducible characters of $K\Sigma'(W_n)$ is given 
by \ref{valeur caractere}. 

\bigskip

\soussection{Radical} The aim of this subsection is to describe the radical 
of $K\Sigma'(W_n)$ in full generality. If $p=0$, then this is done 
in Theorem B (d). Let $\Comp_p(n)=\{C \in \Comp(n)~|~p$ divides $|\WC(C)|\}$. The 
next result must be compared with \cite[Theorem 3]{apv}:

\begin{theo}\label{radical p}
If $K$ is a field of characteristic $p$, then 
$$\Rad K\Sigma'(W_n)=\Ker \th_n^K + \sum_{C \in \Comp_p(n)} K x_C.$$
\end{theo}

\begin{proof}
Let $\IC=\Ker \th_n^K + \sum_{C \in \Comp_p(n)} K x_C$. By Proposition \ref{simples}, 
we get that 
$$\Rad K\Sigma'(W_n)=\bigcap_{\l \in \Bip(n)} \Ker \pi_\l^K.$$
Now, if $\l \in \Bip(n)$, then $\Ker \th_n^K \subset \Ker \pi_\l^K$ and 
$x_C \subset \Ker \pi_\l^K$ for every $C \in \Comp_p(n)$ by \ref{valeur caractere} 
and by Lemma \ref{divisible} (e). Therefore, $\IC \subset \Rad K\Sigma'(W_n)$. 

Now, let $x \in \Rad K\Sigma'(W_n)$. We want to prove that 
$x \in I$. Let $C \in \Comp(n)$ be maximal (for the preorder $\subset_\lamb$) 
such that $\xi_C(x) \neq 0$. By an easy induction argument (on the preorder 
$\subset_\lamb$) we only need to prove that $x' = \sum_{C' \equiv C} \xi_{C'}(x) x_{C'}$ 
belongs to $\IC$. Let $\l=\lamb(C)$. Then, by Lemma \ref{divisible} (d) and 
by \ref{valeur caractere}, we have 
$$0=\pi_\l(x)=\pi_\l(x')=|\WC(C)|\sum_{C' \equiv C} \xi_{C'}(x).$$
To cases may occur:

\medskip

$\bullet$ If $C \in \Comp_p(n)$, then $x' \in \sum_{D \in \Comp_p(n)} K x_D \subset \IC$.

\medskip

$\bullet$ If $C \not\in \Comp_p(n)$, then $p$ does not divide $|\WC(C)|$, 
so $\sum_{C' \equiv C} \xi_{C'}(x) = 0$. So $x' \in \Ker \th_n^K \subset \IC$. 
This completes the proof of the proposition.
\end{proof}

\begin{coro}\label{inclusion}
$|\Bip(n)|=|\Bip_{p'}(n)|+|\lamb(\Comp_p(n))|$.
\end{coro}

\begin{proof}
By Proposition \ref{radical p}, we get that 
$$\dim_K\bigl(\Rad K\Sigma'(W_n)\bigr)=\dim_K(\Ker \th_n^K) + |\lamb(\Comp_p(n))|.$$
On the other hand, we have 
$$\dim_K K\Sigma'(W_n) = \dim_K (\Ker \th_n^K) + |\Bip(n)|$$
$$\dim_K K\Sigma'(W_n) = \dim_K\bigl(\Rad K\Sigma'(W_n)\bigr) + |\Bip_{p'}(n)|.
\leqno{\text{and}}$$
(the last equality follows from Proposition \ref{simples} (b)). 
The corollary now follows from these observations. 
\end{proof}

\bigskip

Note that the above corollary could have been proved directly by a pure 
combinatorial argument. Let us sketch it here. First, 
a bipartition $\l=(\l^+,\l^-)$ is said {\it $p$-regular} 
(respectively {\it $p$-singular}) if it does not belong 
(respectively if it belongs) to $\lamb(\Comp_p(n))$. The set of 
$p$-regular partitions of $n$ (which will be denoted by $\Bip_{p-\reg}(n)$) 
can be described more concretely as follows. 
If $i \ge 1$, we denote by $r_i^+(\l)$ (respectively $r_i^-(\l)$) the number 
of occurrences of $i$ as a part of $\l^+$ (respectively $\l^-$). 
Similarly, if $C \in \Comp(n)$, we denote by $r_i^+(C)$ (respectively $r_i^-(C)$) 
the number of occurrences of $i$ (respectively $-i$) as a part of $C$. 
In other words, $r_i^+(C)=r_i^+(\lamb(C))$ and $r_i^-(C)=r_i^-(\lamb(C))$.
It is readily seen that
\equat\label{normalisateur}
\WC(C) \simeq N_{W_n}(W_C)/W_C \simeq \SG_{r_1^+(C)} \times \dots \times \SG_{r_n^+(C)} 
\times W_{r_1^-(C)} \times \dots \times W_{r_n^-(C)}.
\endequat
Consequently, 
\equat\label{bip 2}
\Bip_{2-\reg}(n)=\{\l \in \Bip(n)~|~\forall i \ge 1,~
r_i^+(\l) \le 1 \text{ and } r_i^-(\l)=0\}
\endequat
and, if $p$ is an odd prime number, 
\equat\label{bip p}
\Bip_{p-\reg}(n)=\{\l \in \Bip(n)~|~\forall i \ge 1,~
r_i^+(\l) \le p-1 \text{ and } r_i^-(\l) \le p-1\}.
\endequat
Now, recall from \S\ref{section conjugacy} that the order $o(\l)$ of $\cox_\l$ 
is equal to the lowest common multiple of $(2\l_1^+,\dots,2\l_r^+,\l_1^-,\dots,\l_s^-)$, 
where $\l^+=(\l_1^+,\dots,\l_r^+)$ and $(\l_1^-,\dots,\l_s^-)$. Therefore, 
\equat\label{bip 2'}
\Bip_{2'}(n)=\{\l \in \Bip(n)~|~\forall i \ge 1,~r_i^+(\l)=r_{2i}^-(\l)=0\}
\endequat
and, if $p$ is an odd prime number, 
\equat\label{bip p'}
\Bip_{p'}(n)=\{\l \in \Bip(n)~|~\forall i \ge 1,~r_{pi}^+(\l)= r_{pi}^-(\l)=0\}.
\endequat
Now, Corollary \ref{inclusion} asserts that 
\equat\label{egalite cardinaux}
|\Bip_{p-\reg}(n)|=|\Bip_{p'}(n)|.
\endequat
This can be proved directly by using the descriptions \ref{bip 2}, \ref{bip p}, \ref{bip 2'} 
and \ref{bip p'} of both sets and by using the classical argument for 
proving the analogue of \ref{egalite cardinaux} for partitions instead 
of bipartitions.

\bigskip

\soussection{Character table}
Let us now talk about the character table of $K\Sigma'(W_n)$. 
By Theorem \ref{radical p}, the classes of the elements of the family 
$(x_{\lamh})_{\l \in \Bip_{p-\reg}(n)}$ 
in the semisimple quotient 
$K\Sigma'(W_n)/\Rad(K\Sigma'(W_n))$ form a $K$-basis of this last space 
(recall that $\lamh$ has been defined in \S\ref{section signed}: it is a representative 
of $\lamb^{-1}(\l)$). Therefore, to compute an irreducible character 
of $K\Sigma'(W_n)$, we only need to give the values on 
$(x_{\lamh})_{\l \in \Bip_{p-\reg}(n)}$. 
We call the {\it character table} of $K\Sigma'(W_n)$ 
the square matrix 
$(\pi_\l^K(x_{\muh}))_{\l \in \Bip_{p'}(n), \mu \in \Bip_{p-\reg}(n)}$. 
By \ref{valeur caractere}, we have~:
\equat\label{triangular}
\text{\it The character table of $K\Sigma'(W_n)$ is upper triangular}
\endequat
(for the order $\subset$ on $\Bip(n)$). 

The character tables of $\QM\Sigma'(W_2)$ and $\QM\Sigma'(W_3)$ will be given at 
the end of this paper. If $p > 0$, the character table of $K\Sigma'(W_2)$ and 
$K\Sigma'(W_3)$ are obtained from the previous ones by reduction modulo $p$ 
and by deleting the appropriate rows and columns.

\bigskip

\section{Loewy length}

\medskip

Recall that the Loewy length of a finite dimensional $K$-algebra $A$ 
is the smallest natural number $r \ge 1$ such that 
$(\Rad A)^r=0$. 
In this section, we shall use the description of the radical 
obtained in Theorem \ref{radical p} to compute the Loewy length 
of $K\Sigma'(W_n)$ (except if $p=2$). But before doing this, we determine the 
Loewy length of the algebra $K \Irr W_n$:

\begin{prop}\label{loewy groth}
The Loewy length of $K\Irr W_n$ is equal to $\begin{cases} 1, & \text{if $p =0$;} \\
n+1, &\text{if $p = 2$;}\\
[n/p]+1, & \text{if $p>2$.}\\
\end{cases}$
\end{prop}

\begin{proof}
The result is obvious if $p=0$ so we may, and we will, assume that 
$p > 0$. If $G$ is a finite group, we denote by $\ell_p(G,1)$ the Loewy length 
of the principal block of $K\Irr G$ (see \cite[\S 3]{bonnafe gr} for the definition 
of the principal block of $K\Irr G$: it is the unique block 
on which the degree map $\deg : K \Irr G \to K$ is non-zero). 
We denote by $\ell_p(n)$ the Loewy length of $K\Irr W_n$. Then 
$$\ell_p(n) \ge \ell_p(W_n,1).\leqno{(1)}$$

If $p$ is odd, then it follows from \cite[Proposition 4.7 (d)]{bonnafe gr} 
that $\ell_p(W_n,1)=\ell_p(\SG_n,1)$. But, by \cite[Corollary 5.8]{bonnafe gr}, 
we have $\ell_p(\SG_n,1)=[n/p]+1$. 
On the other hand, since $W_n$ is isomorphic to a subgroup of $\SG_{2n}$ 
of odd index, it follows from \cite[Proposition 4.7 (a)]{bonnafe gr} 
that $\ell_2(W_n,1) \ge \ell_2(\SG_{2n},1)$. But, by \cite[Corollary 5.8]{bonnafe gr}, 
we have $\ell_2(\SG_{2n},1)=n+1$. So, by using (1), we have proved 
that 
$$\ell_p(n) \ge \begin{cases} 
n+1, &\text{if $p = 2$;}\\
[n/p]+1, & \text{if $p>2$.}\\
\end{cases}\leqno{(2)}$$

We shall now prove that these inequalities are actually equalities. 
For this, we shall need some notation. 
Recall that the $p$-rank of a finite group $G$ is the maximal 
possible rank of an elementary abelian $p$-subgroup of $G$. It will 
be denoted by $\rank_p(G)$. We have
$$\rank_p(\SG_n)=[n/p]\qquad\text{and}\qquad \rank_p(W_n)=\begin{cases}
n & \text{if $p=2$},\\
[n/p] & \text{if $p > 2$.}
\end{cases}$$
If $\l \in \Bip(n)$, we set 
$$\ph_\l=\th_n^K(x_\lamh)\qquad\text{and}\qquad \rank_p(\l)=\rank_p(\WC(\lamh)).$$
In other words, by \ref{normalisateur}, we have
$$\rank_p(\l)=\begin{cases}
\DS{\sum_{i \ge 1} \bigl([r_i^+(\l)/2] + r_i^-(\l)\bigr)} & \text{if $p = 2$,}\\
&\\
\DS{\sum_{i \ge 1} \bigl([r_i^+(\l)/p] + [r_i^-(\l)/p]\bigr)} & \text{if $p > 2$.}
\end{cases}$$
In particular, $\pi_2(\l) \in \{0,1,2,\dots n\}$ and, if $p$ is odd, 
then $\pi_p(\l) \in \{0,1,2,\dots,[n/p]\}$. Note that $\l \in \Bip_{p-\reg}(n)$ 
if and only if $\pi_p(\l)=0$. Note also that $(\ph_\l)_{\l \in \Bip(n)}$ is a 
$K$-basis of $K\Irr W_n$ (see Theorem B). 
Now, by (2), it is sufficient to show that, if $i \ge 0$, then 
$$\bigl(\Rad(K\Irr W_n)\bigr)^i \subset \mathop{\oplus}_{\rank_p(\l) \ge i} K \ph_\l.
\leqno{(3)}$$
So let us now prove (3). Let $\IC_i=\DS{\mathop{\oplus}_{\rank_p(\l) \ge i}} K \chi_\l$. 
We denote by $\IC_i \IC_j$ the space of $K$-linear combinations of elements 
of the form $xy$, where $x \in \IC_i$ and $y \in \IC_j$. Then
$$\IC_i\IC_j \subset \IC_{i+j}.\leqno{(4)}$$
\begin{quotation}
\begin{proof}[Proof of (4)]
We proceed as in \cite[proof of $(\clubsuit)$]{bonnafe sym}. 
For simplification, we set $\NC_C=N_{W_n}(W_C)$ for every $C \in \Comp(n)$. 
We have, for $\l$, $\mu \in \Bip(n)$, 
$$\ph_\l\ph_\mu=\sum_{d \in X_{\lamh\muh}} \ph_{\lamb(\lamh \cap \lexp{d}{\muh})}.$$
The group $\WC(\lamh) \times \WC(\muh)$ acts on $X_{\lamh\muh}$ 
($\WC(\lamh)$ acts by left multiplication while $\WC(\muh)$ 
acts by right multiplication). If $d \in X_{\lamh\muh}$ and 
$(x,y) \in \WC(\lamh)\times \WC(\muh)$, then 
the groups $W_\lamh \cap \lexp{d}{W_\muh}$ and 
$W_\lamh \cap \lexp{xdy^{-1}}{W_\muh}=\lexp{x}{(W_\lamh \cap \lexp{d}{W_\muh})}$ 
are conjugate. In other words, 
$$\lamb(\lamh \cap \lexp{d}{\muh})=\lamb(\lamh \cap \lexp{xdy^{-1}}{\muh}).$$
So, if $X_{\lamh\muh}'$ denotes a set of representatives of 
$(\WC(\lamh)\times\WC(\muh))$-orbits in $X_{\lamh\muh}$, then 
$$\ph_\l\ph_\mu=\sum_{d \in X_{\lamh\muh}'} n_{\l,\mu,d} 
\ph_{\lamb(\lamh \cap \lexp{d}{\muh})},$$
where $n_{\l,\mu,d}$ denotes the cardinality of the orbit of $d$. 
So it is sufficient to show that, if $p$ does not divide $n_{\l,\mu,d}$, then 
$\rank_p(\lamb(\lamh \cap \lexp{d}{\muh})) \ge \rank_p(\l) + \rank_p(\mu)$.

So let $d \in X_{\lamh\muh}'$ be such that $p$ does not divide $n_{\l,\mu,d}$. Let 
$$\fonction{\D_d}{\NC_\lamh \cap \lexp{d}{\NC_\muh}}{\NC_\lamh 
\times \NC_\muh}{w}{(w,d^{-1}wd).}$$
Let $\Delt_d : \NC_\lamh \cap \lexp{d}{\NC_\muh} \to \WC(\lamh) \times 
\WC(\muh)$ denote the composition of $\D_d$ with the canonical 
projection. Then $\Delt_d$ induces an injective morphism 
$\bar{\D}_d : \WC(\l,\mu,d) \to \WC(\lamh) \times \WC(\mu)$, where 
$\WC(\l,\mu,d) = (\NC_\lamh \cap \lexp{d}{\NC_\muh})/W_{\lamh \cap \lexp{d}{\muh}}$. 
Then it is easily checked that $\bar{\D}_d(\WC(\l,\mu,d))$ is the stabilizer 
of $d$ in $\WC(\lamh) \times \WC(\muh)$. In particular, 
$$n_{\l,\mu,d} = \frac{|\WC(\lamh)|.|\WC(\muh)|}{|\WC(\l,\mu,d)|}.$$
So, since $p$ does not divide $n_{\l,\mu,d}$, this means 
that, if $P$ a Sylow $p$-subgroup of $\WC(\l,\mu,d)$, then $\bar{\D}_d(P)$ 
is a Sylow $p$-subgroup of $\WC(\lamh) \times \WC(\muh)$. 
In particular, $\rank_p(\WC(\l,\mu,d)) \ge \rank_p(\l) + \rank_p(\mu)$. 
Now, $\WC(\l,\mu,d)$ is a subgroup of $\WC(\lamh \cap \lexp{d}{\muh})$. 
So $\rank_p(\lamb(\lamh \cap \lexp{d}{\muh})) \ge \rank_p(\l)+\rank_p(\mu)$, 
as desired.
\end{proof}
\end{quotation}

\medskip

By (4), $\IC_i$ is an ideal of $K\Irr W_n=\IC_0$. Moreover, 
again by (4), $\IC_1$ consists of nilpotent elements. So 
$\IC_1 \subset \Rad(K\Irr W_n)$. 
On the other hand, $\dim_K \IC_1 = |\Bip(n)|-|\Bip_{p-\reg}(n)|$, 
so $\dim_K \IC_1 = \dim_K\bigl(\Rad(K\Irr W_n)\bigr)$ 
(see \cite[Corollary 2.16]{bonnafe gr}). 
So $\IC_1=\Rad(K\Irr W_n)$. But, by (4), $\IC_1^i \subset \IC_i$. 
This shows (3), so the proof of the proposition is complete.
\end{proof}

We are now ready to prove the main theorem of this section 
(compare with \cite[\S 5.E]{BP}):

\begin{theo}\label{loewy}
If $p \neq 2$, then 
the Loewy length of $K\Sigma'(W_n)$ is $n$. If $p=2$, then 
this Loewy length lies in $\{n,n+1,\dots,2n-1\}$.
\end{theo}

\begin{proof} 
Let $l_p(n)$ denote the Loewy length of $K \Sigma'(W_n)$. 
If $n = 1$, then the result of the Theorem is easily checked. So 
we may, and we will, assume that $n \ge 2$. 
The proof will proceed in two steps. 

\medskip

$\bullet$ {\it First step: upper bound.} 
We use here the notation of Example \ref{exemple bilateres}. 
Let us first prove the following result:
if $k \ge 0$ and if $x \in \Rad K\Sigma'(W_n)$, then:

$$x K\Sigma_{\FC_k^-(n)}'(W_n) \subset K\Sigma_{\FC_{k+1}^-(n)}'(W_n);\leqno{(1)}$$

$$\text{\it If $p \neq 2$, then 
$x K\Sigma_{\FC_k(n)}'(W_n) \subset K\Sigma_{\FC_{k+1}(n)}'(W_n)$.}\leqno{(2)}$$

\medskip

\begin{quotation}
\begin{proof}[Proof of (1) and (2)] 
Let $A_n$ denote the algebra 
$K\Sigma'(W_n)/K\Sigma_{\FC_1^-(n)}'(W_n)$. Recall that 
$K\Sigma_{\FC_1^-(n)}'(W_n)$ 
is a two-sided ideal of $K\Sigma'(W_n)$ (see Example \ref{exemple bilateres}). 
Then $A_n \simeq K[T]/(T(T-2))$, where $T$ 
is an indeterminate. Indeed, $A_n$ has dimension $2$ and is generated 
by the image $t_n$ of $x_{(-n)}$ and it is easily checked that $t_n^2=2t_n$ 
(this follows for instance from Proposition C (c), from Proposition 
\ref{tau bis}, from Lemma \ref{divisible} (c) and from the fact that 
$|N_{W_n}(\SG_n)/\SG_n|=|\WC(-n)|=2$). In particular, if 
$p \neq 2$, then $A_n \simeq K \times K$ is split semisimple. 

Now, let $D \in \Comp(n)$. 
Write $D=(d_1,\dots,d_r)$ and let $a=\Res_{D^+}(x)$. Then 
$$xx_D=x_{D^+} a x_D^{D^+}.$$
Since $\Res_{D^+}$ is a morphism of algebras, $a$ is 
a nilpotent element of the algebra $\Sigma'(W_{D^+})\simeq 
\Sigma'(W_{|d_1|}) \otimes \dots \otimes \Sigma'(W_{|d_r|})$. 
In particular, its image $\aba$ in $A_{D^+}=A_{|d_1|} \otimes \dots \otimes A_{|d_r|}$ 
is also nilpotent. So, if $D \in \FC_k^-(n)$ (respectively if $D \in \FC_k(n)$ and 
$p \neq 2$) then the above description of $A_n$ shows that 
$xx_D \in K\Sigma_{\FC_{k+1}^-(n)}'(W_n)$ (respectively 
$xx_D \in K\Sigma_{\FC_{k+1}(n)}'(W_n)$). 
\end{proof}
\end{quotation}

\medskip

Since $\FC_n(n)=\FC_{2n}^-(n)=\vide$, then the statement (1) above shows that 
$l_2(n) \le 2n$ and the statement (2) shows that, if $p \neq 2$, then $l_p(n) \le n$. 
We shall show now that, if $n \ge 2$, then $l_2(n) \le 2n-1$. 
So let $a_1$,\dots, $a_{2n-1} \in \Rad K\Sigma'(W_n)$. Then, by (a), we have 
$$a_1 \dots a_{2n-1} \in K\Sigma_{\FC_{2n}(n)}'(W_n)= K x_{(-1,-1,\dots,-1)}.$$
Let $\l \in K$ be such that $a_1 \dots a_{2n-1}=\l x_{(-1,-1,\dots,-1)}$. 
Then $\th_n^K(a_1 \dots a_{2n-1}) = \l \chi_n$, where $\chi_n$ is the 
regular character of $W_n$. But, since $\th_n^K(a_i)$ belongs to the radical 
of the $K$-algebra $K \Irr W_n$, since this algebra has Loewy length 
$\le n+1$ (see Proposition \ref{loewy groth}) and since $n+1 \le 2n-1$ 
(because $n \ge 2$), we get that $\l =0$, as desired. 
So we have proved the following results:

$$\text{\it If $n \ge 2$, then $l_2(n) \le 2n-1$.}\leqno{(3)}$$

$$\text{\it If $p \neq 2$, then $l_p(n) \le n$.}\leqno{(4)}$$

\medskip

$\bullet$ {\it Second step: lower bound.} 
Let $a=x_{(n-1,-1)}-x_{(-1,n-1)}$. Then $a \in \Ker \th_n^K$. 
If $1 \le i \le n$, let $C_i$ denote the signed compositions 
$(1,\dots,1,-1,1,\dots,1)$ of $n$, where the $-1$ term is in position $i$. Then:

$$a^{n-1} = 
\sum_{i=1}^n (-1)^i \begin{pmatrix} n-1 \\ i-1 \end{pmatrix} x_{C_i}.\leqno{(5)}$$

\begin{quotation}
\begin{proof}[Proof of (5)]
If $1 \le j \le n$ and if $1 \le i \le n+1-j$, we denote 
by $C_{i,j}$ the signed composition $(-1,\dots,-1,j,-1,\dots,-1)$, 
where $j$ appears in the $i$-th position (for instance, $C_{i,1}=C_i$). 
For simplification, we set $s_{i,j}=x_{C_{i,j}}$. 
We have in particular $a=s_{1,n-1}-s_{2,n-1}$. We want to show 
by induction on $k \in \{1,2,\dots,n-1\}$ that 
$$a^k = \sum_{i=1}^{k+1} (-1)^i \begin{pmatrix} k \\ i-1 \end{pmatrix} s_{i,n-k}.
\leqno{(5^+)}$$
Note that the formula (5) is obtained by specialising $k$ to $n-1$ in the 
formula $(5^+)$. For proving $(5^+)$ by induction, it is sufficient to show 
that
$$s_{1,n-1} s_{i,j} = \a_{i,j} s_{i,j} + s_{i,j-1} 
\quad\text{and}\quad
s_{2,n-1} s_{i,j}  = \a_{i,j} s_{i,j} + s_{i+1,j-1}$$
for some $\a_{i,j} \in \NM$. The first equality is easily 
checked using the description of $X_{n-1,-1}$ given in the proof 
of Proposition \ref{surjectif B}. The second one follows from a 
similar argument. 
\end{proof}
\end{quotation}

\medskip

The statement (5) above shows that $l_p(n) \ge n$. By (3) and (4), the proof 
of the Theorem is complete.
\end{proof}

\remark{l2}
Keep the notation of the proof of the previous Theorem. 
It is probable that $l_2(n)=2n-1$ whenever $n \ge 2$ 
(note that $l_2(1)=2$, $l_2(2)=3$, $l_2(3)=5$, $l_2(4)=7$ and $l_2(5)=9$). 
In fact, it is probable 
that the element $a$ defined in the above proof lies in 
$(\Rad \FM_{\! 2}\Sigma'(W_n))^2$ 
(it has been checked for $n \le 5$): this would imply that $l_2(n)=2n-1$ for $n \ge 2$ 
(see the statement (5) of the above proof).\finl

%
%
%
%
%
%

\bigskip

\section{Projective modules, Cartan matrix}

\medskip

\soussection{Projective modules} 
If $\l \in \Bip(n)$, we denote by $e_\l^\QM : W_n \to \QM$ the characteristic 
function of $\CC(\l)$. We may, and we will, view it as an element of 
$\QM\Irr W_n$: we have
\equat
e_\l^\QM = \frac{|\CC(\l)|}{|W_n|}\sum_{\chi \in \Irr W_n} \chi(\cox_\l) \chi.
\endequat
Then $(e_\l^\QM)_{\l \in \Bip(n)}$ is a family of orthogonal primitive idempotents 
of $\QM\Irr W_n$ such that $\sum_{\l \in \Bip(n)} e_\l^\QM = 1_n$. 
Since the morphism $\th_n^\QM$ is surjective, there exists \cite[Theorem 3.1 (f)]{thevenaz} 
a family $(E_\l^\QM)_{\l \in \Bip(n)}$ of primitive idempotents of $\QM\Sigma'(W_n)$ such that
\begin{quotation}
\begin{itemize}
\itemth{1} $\forall~\l \in \Bip(n)$, $\th_n^\QM(E_\l^\QM)=e_\l^\QM$.

\itemth{2} $\forall~\l,\mu \in \Bip(n)$, $E_\l^\QM E_\mu^\QM = E_\mu^\QM E_\l^\QM = 
\d_{\l\mu} E_\l^\QM$.

\itemth{3} $\sum_{\l \in \Bip(n)} E_\l^\QM = 1$.
\end{itemize}
\end{quotation}
Let $\PC_\l^\QM=\QM\Sigma'(W_n) E_\l^\QM$. It is an indecomposable projective 
$\QM\Sigma'(W_n)$-module: this is the projective cover of $\DC_\l^\QM$. 
Moreover, 
\equat
\mathop{\oplus}_{\l \in \Bip(n)} \PC_\l^\QM = \QM\Sigma'(W_n).
\endequat
If $p=0$, then $\QM \subset K$ and we set $E_\l^K=E_\l^\QM$. 
Note that, if $p=0$, then $(E_\l^K)_{\l \in \Bip(n)}$ is still a family 
of orthogonal primitive idempotents of $K\Sigma'(W_n)$ (since $\QM\Sigma'(W_n)$ 
is split). 

\bigskip

\example{x un un prime idempotent}
We keep the notation introduced in Example \ref{x un un prime}. 
If $I \subset I_n^+$, then 
$$\lamb(\Cb(I))=((\underbrace{1,1,\dots,1}_{\text{$|I|$ times}}), 
(\underbrace{1,1,\dots,1}_{\text{$n-|I|$ times}})).$$
Then
\equat\label{conjugue x un un}
\text{\it the idempotent $x_{\Cb(I)}'/(2^{n-|I|}|\SG_n(I)|)$ 
is conjugate to $E_{\lamb(\Cb(I))}^\QM$.}
\endequat
Let us prove this result.
Since $x_{\Cb(I)}'/(2^{n-|I|}|\SG_n(I)|)$ is an idempotent (see \ref{x un un quasi}), 
it is sufficient to show that 
$\th_n^\QM(x_{\Cb(I)}'/(2^{n-|I|}|\SG_n(I)|))=\th_n^\QM(E_{\lamb(\Cb(I))}^\QM)$ 
(see \cite[Theorem 3.1 (e)]{thevenaz}). In other words, it is sufficient 
to show that $\th_n^\QM (x_{\Cb(I)}')$ is a multiple of $e_{\lamb(\Cb(I)}^\QM$. 

If $J \subset I_n^+$, let $\TG_J$ denote the subgroup of $\TG_n$ 
generated by $(t_i)_{i \in J}$ and let $t_J=\prod_{i \in J} t_i$. 
Then $t_J \in \CC(\lamb(\Cb(J)))$ and  
$$\th_n^\QM(x_{\Cb(I)}')=\Ind_{\TG_I}^{W_n} f_I,\quad
\text{where}\quad 
f_I=\sum_{J \subset I} \Bigl(-\frac{1}{2}\Bigr)^{|I|-|J|} 1_J.$$
Here, $1_J$ denotes the trivial character of $\TG_J$. It is now 
sufficient to show that $f_I$ is the characteristic function of $\{t_I\}$ 
in $\TG_I$: since $\TG_I$ is an elementary abelian $2$-group, this is easily 
reduced, by direct products, to the case where $|I|=1$ (for which it is obvious).\finl

\bigskip

Let us now assume that $p > 0$. 
For each $\l \in \Bip_{p'}(n)$, we denote by $\CC_{p'}(\l)$ the set of elements 
$w$ in $W_n$ such that $w_{p'}$ belongs to $\CC(\l)$. It is a union of conjugacy classes 
of $W_n$. We set
$$e_{\l,p'}^\QM =\sum_{\substack{\mu \in \Bip(n) \\ \mu_{p'} = \mu}} e_\mu^\QM.$$
This is the characteristic function of $\CC_{p'}(\l)$. It is an idempotent 
of $\QM \Irr W_n$ and, by \cite[Corollary 2.21]{bonnafe gr}, it is a primitive 
idempotent of $\ZM_{(p)} \Irr W_n$. We denote by $e_\l^K$ its image 
in $K\Irr W_n$: it is still a primitive idempotent of $K\Irr W_n$. 
Then $(e_\l^K)_{\l \in \Bip_{p'}(n)}$ is a family of orthogonal primitive idempotents 
of $K\Irr W_n$ such that $\sum_{\l \in \Bip_{p'}(n)} e_\l^K = 1_n$. 
Since the morphism $\th_n^K$ is surjective, there exists 
\cite[Theorem 3.1 (f)]{thevenaz} a family 
$(E_\l^K)_{\l \in \Bip_{p'}(n)}$ of primitive idempotents of $K\Sigma'(W_n)$ such that
\begin{quotation}
\begin{itemize}
\itemth{1} $\forall~\l \in \Bip_{p'}(n)$, $\th_n^K(E_\l^K)=e_\l^K$.

\itemth{2} $\forall~\l,\mu \in \Bip_{p'}(n)$, $E_\l^K E_\mu^K = E_\mu^K E_\l^K = 
\d_{\l\mu} E_\l^K$.

\itemth{3} $\sum_{\l \in \Bip(n)} E_\l^K = 1$.
\end{itemize}
\end{quotation}
Let $\PC_\l^K=K\Sigma'(W_n) E_\l^K$. It is an indecomposable projective 
$K\Sigma'(W_n)$-module: this is the projective cover of $\DC_\l^K$. 
Moreover, 
\equat
\mathop{\oplus}_{\l \in \Bip_{p'}(n)} \PC_\l^K = K\Sigma'(W_n).
\endequat

We conclude this section by a useful remark on the idempotents 
$E_\l^\QM$: this will be used for proving the unitriangularity of the 
Cartan matrix of $\QM\Sigma'(W_n)$.

\begin{prop}\label{localisation}
Let $D \in \Comp(n)$, let $\l = \lamb(D)$ and let $D'$ be a semi-positive 
signed composition of $n$ such that $D \subset D'$. Then there exists a primitive 
idempotent $E$ of $\QM\Sigma'(W_n)$ satisfying the following two conditions:
\begin{itemize}
\itemth{a} $\th_n^\QM(E)=e_\l^\QM$.

\itemth{b} $E \in \QM\Sigma'(W_n) x_{D'}$.
\end{itemize}
In particular, $E$ is conjugate to $E_\l^\QM$. 
\end{prop}

\remark{choix D}
In the above Proposition, one can choose $D'=D^+$.\finl

\begin{proof}
For simplification, let $\KC=\Ker(\Res_{D'})$ and $\IC=\QM\Sigma'(W_n) x_{D'}$. 
By Proposition \ref{noyau} (b), we have 
$$\QM\Sigma'(W_n)=\KC \oplus \IC.$$
In particular, the restriction of $\Res_{D'}$ to the left ideal 
$\IC$ is injective. Moreover, as a direct consequence 
of the hypothesis, we get that $\Res_{W_{D'}}^{W_n} e_\l^\QM \neq 0$, 
so in particular $\Res_{D'} E_\l^\QM \neq 0$ (see also Proposition \ref{restriction} (d)). 
Let us write $E_\l^\QM = A + E$, with $A \in \KC$ and 
$E \in \IC$. Then $E^2-E \in \IC$ and 
$\Res_{D'}(E^2-E)=\Res_{D'}((E_\l^\QM)^2-E_\l^\QM)=0$. Therefore, 
$E^2=E$. Moreover, $AE \in \KC \cap \IC$, so $AE=0$. In other words, 
$E_\l^\QM E=E^2=E$. This shows in particular that 
$$\dim_\QM \QM\Sigma'(W_n) E = \dim_\QM \QM\Sigma'(W_n) E_\l^\QM E \le 
\dim_\QM \QM\Sigma'(W_n) E_\l^\QM.\leqno{(*)}$$
Now, $\Res_{D'}(E_\l^\QM)=\Res_{D'}(E)$. Since $E_\l^\QM$ is primitive, this implies 
that $E=\EC_\l+F$, where $\EC_\l$ and $F$ are orthogonal idempotent and 
$\EC_\l$ is conjugate to $E_\l^\QM$ (see \cite[Theorem 3.2 (c)]{thevenaz}). 
But, by $(*)$, we get that $F=0$, so that $\th_n^\QM(E)=e_\l^\QM$. 
This shows the proposition. 
\end{proof}

\soussection{About the structure of ${\boldsymbol{K W_n}}$ as a 
left ${\boldsymbol{K \Sigma'(W_n)}}$}
The next result is the analogue of \cite[Theorem 7.15]{bbht}:

\begin{prop}\label{dimension proj}
If $p=0$ and if $\l \in \Bip(n)$, then $\dim_K K W_n E_\l^K = |\CC(\l)|$.
\end{prop}

\begin{proof}
We may assume that $K=\QM$. 
Let $\TC_n : \QM W_n \to \QM$ denote the unique linear map such that 
$\TC_n(1)=1$ and $\TC_n(w)=0$ for every $w \in W_n$ which is different 
from $1$. Then $\TC_n$ is the canonical symmetrizing form on $\QM W_n$. 
Now, if $x \in \QM W_n$, then the trace of the multiplication 
by $x$ on $\QM W_n$ (on the left or on the right) is equal to $|W_n|\TC_n(x)$. 
Therefore, since $E_\l^\QM$ is an idempotent, $\dim_\QM \QM W_n E_\l^\QM=|W_n| 
\TC_n(E_\l^\QM)$. 
But, by \cite[Proposition 3.8]{BH}, we have 
$$\TC_n(E_\l^\QM) = \langle \th_n^\QM(E_\l^\QM),\th_n^\QM(1) \rangle_{W_n} = 
\langle e_\l^\QM,1_n \rangle_{W_n} = \frac{|\CC(\l)|}{|W_n|},$$
as expected.
\end{proof}

In the same spirit, we have the following result:

\begin{prop}\label{caractere reguliere}
The character of the $\QM\Sigma'(W_n)$-module $\QM W_n$ 
is $\DS{\sum_{\l \in \Bip(n)}} |\CC(\l)| \pi_\l^\QM$.
\end{prop}

\begin{proof}
If $C \in \Comp(n)$, the trace of $x_C$ in its left action on $\QM W_n$ 
is equal to $|W_n|\TC_n(x_C)=|W_n|$. On the other hand, 
\eqna
\DS{\sum_{\l \in \Bip(n)}} |\CC(\l)|\pi_\l^\QM(x_C) &=& 
\DS{\sum_{\l \in \Bip(n)}} |\CC(\l)| \th_n^\QM(x_C)(\cox_\l) \\
&=& \DS{\sum_{\l \in \Bip(n)}} |W_n| \langle \th_n^\QM(x_C),e_\l^\QM \rangle_{W_n} \\
&=& |W_n| \langle \th_n^\QM(x_C),1_n \rangle_{W_n} \\
&=& |W_n|,
\endeqna
as desired.
\end{proof}

\begin{prop}\label{dim p}
If $p > 0$ and if $\l \in \Bip_{p'}(n)$, then $\dim_K KW_n E_\l^K = |\CC_{p'}(\l)|$.
\end{prop}
\def\gfp{{\FM_{\! p}}}

\begin{proof}
We may, and we will, assume that $K=\gfp$. The idempotent $E_\l^\gfp$ can be lifted 
to an idempotent $E_\l^{\ZM_p}$ of $\ZM_p\Sigma'(W_n)$, where $\ZM_p$ denotes the ring 
of $p$-adic integers \cite[Theorem 3.2 (b)]{thevenaz}. It is sufficient to show that 
$$\dim_{\QM_p} \QM_p\Sigma'(W_n) E_\l^{\ZM_p}=|\CC_{p'}(\l)|.\leqno{(?)}$$
Now, $\th_n(E_\l^{\ZM_p})$ is an idempotent of 
$\ZM_p \Irr W_n$ that lifts $e_\l^K$. Therefore, 
$\th_n(E_\l^{\ZM_p})=e_{\l,p'}^\QM$ by the unicity of liftings 
in commutative algebras \cite[Theorem 3.2 (d)]{thevenaz}. 
Therefore, $E_\l^{\ZM_p}$ is conjugate 
to the idempotent $\sum_{\mu \in \Bip(n), \mu_{p'}=\l} E_\mu^\QM$ 
(see \cite[Theorem 3.2 (d)]{thevenaz}). So the result 
follows from Proposition \ref{dimension proj}.
\end{proof}

\bigskip

\soussection{Cartan matrix} 
We return to the general situation, namely we assume that 
$K$ is a field of characteristic $p \ge 0$. We denote by 
$\Cartan(K\Sigma'(W_n))$ the Cartan matrix of $K\Sigma'(W_n)$. 
It is the square matrix $([\PC_\l^K:\DC_\mu^K])_{\l, \mu \in \Bip_{p'}(n)}$, 
where $[\PC_\l^K:\DC_\mu^K]$ denotes the multiplicity of $\DC_\mu^K$ 
as a chief factor in a Jordan-H\"older series of $\PC_\l^K$. Recall that 
\equat\label{dim hom}
[\PC_\l^K:\DC_\mu^K]=\dim_K \Hom_{K\Sigma'(W_n)}(\PC_\mu^K,\PC_\l^K)
\endequat
and that we have a canonical isomorphism of vector spaces 
\equat\label{hom ef}
\Hom_{K\Sigma'(W_n)}(\PC_\mu^K,\PC_\l^K) \simeq E_\mu^K K\Sigma'(W_n) E_\l^K.
\endequat
Moreover, the isomorphism \ref{hom ef} is an isomorphism of algebras 
whenever $\l=\mu$. 

Let $D_n^K=(\d_{\l_{p'},\mu})_{\l \in \Bip(n),\mu \in \Bip_{p'}(n)}$, where 
$\d_{!,?}$ is the Kronecker symbol. If $p$ does not 
divide the order of $|W_n|$, this is just the identity matrix. 
In general, it may be seen as the decomposition 
matrix from $\QM\Sigma'(W_n)$ to $K\Sigma'(W_n)$ (see \cite[\S 7.4]{geck pfeiffer} 
for the general definition of a decomposition matrix). 
The next lemma reduces the computation of $\Cartan(K\Sigma'(W_n))$ to the computation 
of $\Cartan(\QM\Sigma'(W_n))$ by making use of the decomposition matrix $D_n^K$ 
(see \cite[Theorem 8]{apv} for the analogue of the next result for Solomon 
descent algebras). 

\begin{lem}\label{triangle cde}
We have $\Cartan(K\Sigma'(W_n))=\lexp{t}{D_n^K} \Cartan(\QM\Sigma'(W_n)) D_n^K$.
\end{lem}

\begin{proof}
This follows from \cite[\S 2.3]{GR}.
\end{proof}

The next result is a first decomposition of the Cartan matrix 
of $K\Sigma'(W_n)$ into diagonal blocks (whenever $p \neq 2$), 
according to the action of $w_n$ on simple modules.

\begin{lem}\label{congruence}
Assume that $p \neq 2$. Let $\l$, $\mu \in \Bip_{p'}(n)$. 
If $[\PC_\l^K:\DC_\mu^K]\neq 0$, then 
$$\length^-(\l) \equiv \length^-(\mu) \mod 2.$$
\end{lem}

\begin{proof}
Let $\l$, $\mu \in \Bip_{p'}(n)$ be such that $[\PC_\l^K:\DC_\mu^K]\neq 0$ 
and let $? \in \{+,-\}$. First, note that $e_n^?\PC_\l^K =\PC_\l^K$ if and only if 
$e_n^? \DC_\l^K$ because $\PC_\l^K$ is indecomposable. On the other hand, 
if $e_n^?\PC_\l^K =\PC_\l^K$, then $e_n^? \DC_\mu^K =\DC_\mu^K$. 
So the result follows from \ref{tau signe}.
\end{proof}

The main result of this section is the following:

\begin{theo}\label{cartan triangular}
The Cartan matrix $\Cartan(\QM\Sigma'(W_n))$ is unitriangular. More precisely, 
if $\l$ and $\mu$ are two distinct bipartitions of $n$, then:
\begin{itemize}
\itemth{a} $[\PC_\l^\QM:\DC_\l^\QM]=1$.

\itemth{b} If $[\PC_\l^\QM:\DC_\mu^\QM] \neq 0$, then $\length(\mu) > \length(\l)$.
\end{itemize}
\end{theo}

\begin{proof}
Let $\FC=\{C \in \Comp(n)~|~\length(C) > \length(\l)\}$. Then 
$\FC$ is saturated so $\IC=\QM\Sigma_\FC'(W_n)$ is a two-sided ideal 
of $\QM\Sigma'(W_n)$. The theorem follows from the fact that 
$\PC_\l^\QM \subset \QM E_\l^\QM + \IC$, 
which is an immediate consequence of the statement (2) in the proof 
of the Theorem \ref{loewy}.
\end{proof}

We shall give at the end of this paper the Cartan matrices of 
$\QM\Sigma'(W_2)$, $\QM\Sigma'(W_3)$ and $\QM\Sigma'(W_4)$. 

\begin{coro}\label{centre semisimple}
If $p$ does not divide the order of $W_n$ (i.e. if $p=0$ or $p > \max(2,n)$), 
then the centre of $K\Sigma'(W_n)$ is split semisimple.
\end{coro}

\begin{proof}
Note that $\Bip(n)=\Bip_{p'}(n)$. 
Let $Z$ be the centre of $K\Sigma'(W_n)$. Then the map 
$Z \to \End_K (K\Sigma'(W_n))$ sending $z \in Z$ to the left 
multiplication by $z$ is injective. Moreover, the image is contained 
in $\oplus_{\l \in \Bip_{p'}(n)} \End_{K\Sigma'(W_n)} \PC_\l^K$. 
But, by \ref{dim hom}, by Lemma \ref{triangle cde} (and the fact that the matrix 
$D_n^K$ is the identity) and by Theorem \ref{cartan triangular} (a), we have an 
isomorphism of $K$-algebras 
$$\End_{K\Sigma'(W_n)} \PC_\l^K \simeq K.$$
So $Z$ is a subalgebra of $K \times \cdots \times K$ ($|\Bip(n)|$ times). 
The proof of the corollary is complete. 
\end{proof}

\example{pas semisimple centre}
If $p$ divides the order of $W_n$, then the centre of $K\Sigma'(W_n)$ is not 
semisimple. Indeed, the element $x_{(-1,-1,\dots,-1)}$ is central in 
$\QM\Sigma'(W_n)$ and $(x_{(-1,-1,\dots,-1)})^2=|W_n| x_{(-1,-1,\dots,-1)}$.\finl

\begin{coro}\label{reduction stable}
Let $Z_n^R$ denote the centre of $R\Sigma'(W_n)$. 
If $p$ does not divide the order of $W_n$, then 
the natural map $K \otimes_\ZM Z_n^\ZM \to Z_n^K$ is an 
isomorphism of algebras.
\end{coro}

\begin{proof}
It is sufficient to show that $\dim_K Z_n^K = \dim_\QM Z_n^\QM$. 
But, since $Z_n^K$ is split semisimple (see Corollary \ref{centre semisimple}), 
its dimension is equal to the number of blocks of $K\Sigma'(W_n)$. 
This number is determined by the Cartan matrix of $K\Sigma'(W_n)$. 
Since the Cartan matrices of $K\Sigma'(W_n)$ and $\QM\Sigma'(W_n)$ 
coincide by Lemma \ref{triangle cde}, the result follows. 
\end{proof}

The next example shows that the Corollary \ref{reduction stable} does 
not hold for any $p$ and any $n$. 

\bigskip

\example{dim centre} 
It would be interesting to determine the centre of $K\Sigma'(W_n)$. 
Note that this dimension is always $\ge 4$. Indeed, if $p \neq 2$, 
then $x_n=1$, $w_n$, $x_{\vide}$ and $x_{1,1,\dots,1}'$ are linearly 
independent central elements (see \ref{central x un}). If $p=2$, 
then $x_{1,1,\dots,1}'$ must me replaced by the image of 
$2^{n-1} x_{1,1,\dots,1}' - x_{\vide}/2 \in \ZM\Sigma'(W_n)$ 
in $K\Sigma'(W_n)$. 

The next table, obtained using {\tt CHEVIE} \cite{chevie}, 
provides the dimension of this centre for $n \le 5$: 
it depends on the characteristic $p$ of $K$. 

\def\vertical{\vphantom{\DS{\frac{A}{A}}}}

$$\begin{array}{|c||c|c|c|}
\hline
\vertical n \backslash p & 0 & 2 & \ge 3 \\
\hline
\hline
\vertical 1 & ~2~ & ~2~ & ~2~ \\
\hline
\vertical 2 & 4 & 4 & 4 \\
\hline
\vertical 3 & 4 & 4 & 4 \\
\hline
\vertical 4 & 5 & 6 & 5 \\
\hline
\vertical 5 & 4 & 4 & 4 \\
\hline
\end{array}$$

\medskip

\noindent 
Note that the case $p > n$ has been handled by using the Lemma \ref{reduction stable}. 
It would also be interesting to determine for which pairs $(p,n)$ does the Lemma 
\ref{reduction stable} hold. For instance, does it hold if $p$ is odd?\finl

\bigskip

We conclude this section by proving that the Cartan matrix 
of $\QM\Sigma'(W_n)$ is a submatrix of the Cartan matrix of 
$\QM\Sigma'(W_{n+1})$. We identify $\Bip(n,-1)$ with $\Bip(n)$ 
and the map $\t_{(n,-1)} : \Bip(n) \to \Bip(n+1)$ defined in \S\ref{def restriction} 
will be denoted simply by $\t_n$. Then:

\bigskip

\begin{theo}\label{cartan induction}
Let $\l$ and $\mu$ be two bipartitions of $n$. Then
$$[\PC_{\t_n(\l)}^\QM:\DC_{\t_n(\mu)}^\QM]=[\PC_\l^\QM:\DC_\mu^\QM].$$
\end{theo}

\bigskip

\begin{proof}
By Proposition \ref{localisation}, we may, and we will, assume that 
$E_{\t_n(\l)}^\QM$ and $E_{\t_n(\mu)}^\QM$ belong to 
$\QM\Sigma'(W_{n+1})x_{n,-1}$. In particular, by Proposition 
\ref{noyau} (b), $E_{\t_n(\mu)}^\QM \QM\Sigma'(W_{n+1}) E_{\t_n(\l)}^\QM$ 
is mapped isomorphically to $\EC_\l \QM\Sigma'(W_n) \EC_\mu$ through 
the map $\Res_n^{n+1}$, where $\EC_\l=\Res_n^{n+1}(E_\l^\QM)$ and 
$\EC_\mu=\Res_n^{n+1}(E_\mu^\QM)$. So it remains to show 
that $\EC_\l$ and $\EC_\mu$ are conjugate to $E_\l^\QM$ and $E_\mu^\QM$ 
respectively (see \ref{dim hom} and \ref{hom ef}). 

Let us prove it for $\l$ (this is sufficient). 
By Proposition \ref{restriction} (d), we have 
$$\th_n^\QM(\EC_\l) = \Res_{W_n}^{W_{n+1}} e_{\t_n(\l)}^\QM = e_\l^\QM.$$
But, $\EC_\l$ is a primitive idempotent in the image of $\Res_n^{n+1}$ 
(see \cite[Theorem 3.2 (d)]{thevenaz}) and since $\Res_n^{n+1}$ is surjective 
(see Proposition \ref{surjectif B}), we get that 
$\EC_\l$ is a primitive idempotent of $\QM\Sigma'(W_n)$. 
So $\EC_\l$ and $E_\l^\QM$ are conjugate \cite[Theorem 3.2 (c)]{thevenaz}.
\end{proof} 

\section{Numerical results}

\medskip

For simplification, a bipartition $((\l_1^+,\dots,\l_r^+),(\l_1^-,\dots,\l_s^-))$ 
will be denoted in a compact way 
$\l_1^+\dots\l_r^+;\l_1^-\dots\l_s^-$. For instance, 
$31;411$ stands for $((3,1),(4,1,1))$ and $\vide;221$ stands for 
$((),(2,2,1))$. If $i \ge 1$, the number $-i$ will be denoted 
by $\bar{i}$: for instance, the signed composition $(2,-3,-1,1,-2)$ 
will be denoted by $(2,\bar{3},\bar{1},1,\bar{2})$. 

We shall give here the Cartan matrix and the primitive central idempotents 
of the algebras $\QM\Sigma'(W_n)$ for $n \in \{2,3,4\}$. 
For $n \in \{2,3\}$, we also give the character table and 
an example of a family $(E_\l^\QM)_{\l \in \Bip(n)}$. 
Note that they are obtained by lifting the idempotents 
$(e_\l^\QM)_{\l \in \Bip(n)}$ by using {\tt CHEVIE} \cite{chevie} 
and the algorithm described 
in \cite[Theorem 3.1 (b) and (f)]{thevenaz}.
In the next tables, we have replaced zeroes 
by dots. Note also that, for simplicity, the idempotents will be expressed 
in the basis $(x_C')_{C \in \Comp(n)}$ constructed in Section \ref{action plus long}.

\bigskip

\soussection{The case ${\boldsymbol{n=2}}$} 
The character table of $\QM\Sigma'(W_2)$ is:

$$\begin{array}{|c|ccccc|}
\hline
\vertical  & x_2 & x_{\dbar} & x_{1,1} & x_{1,\ubar} & x_{\ubar,\ubar}  \\
\hline
\vertical \pi_{2;\vide}^\QM       & 1 & . & . & . & . \\
\vertical \pi_{\vide;2}^\QM       & 1 & 2 & . & . & . \\
\vertical \pi_{11;\vide}^\QM      & 1 & . & 2 & . & . \\
\vertical \pi_{1;1}^\QM           & 1 & . & 2 & 2 & . \\
\vertical \pi_{\vide;11}^\QM      & 1 & 4 & 2 & 4 & 8 \\
\hline
\end{array}$$

\bigskip

We can take for the family $(E_\l^\QM)_{\l \in \Bip(2)}$ the following 
idempotents:
\eqna
E_{2;\vide}^\QM & =& \DS{x_2' - \frac{1}{2}x_{1,1}' + \frac{1}{8} x_{\ubar,\ubar}'}\\ 
E_{\vide;2}^\QM &=& \DS{\frac{1}{2} x_{\dbar}' + 
\frac{1}{4}(x_{1,\ubar}' - x_{\ubar,1}')}\\
E_{11;\vide}^\QM &=&\DS{\frac{1}{2}x_{1,1}'}\\
E_{1;1}^\QM &=& \DS{\frac{1}{4}(x_{1,\ubar}'+x_{\ubar,1}') }\\
E_{\vide;11}^\QM&=& \DS{\frac{1}{8} x_{\ubar,\ubar}'}.
\endeqna

\bigskip

The Cartan matrix of $\QM\Sigma'(W_2)$ is given by:

$$\begin{array}{|c|ccccc|}
\hline
\vertical  & \DC_{2;\vide}^\QM & 
\DC_{\vide;2}^\QM & \DC_{1;1}^\QM & \DC_{11;\vide}^\QM & \DC_{\vide;11}^\QM   \\
\hline
\vertical \PC_{2;\vide}^\QM       & 1 & . & . & . & . \\
\vertical \PC_{\vide;2}^\QM       & . & 1 & 1 & . & . \\
\vertical \PC_{1;1}^\QM           & . & . & 1 & . & . \\
\vertical \PC_{11;\vide}^\QM      & . & . & . & 1 & . \\
\vertical \PC_{\vide;11}^\QM      & . & . & . & . & 1 \\
\hline
\end{array}$$

\bigskip

The primitive central idempotents of $\QM\Sigma'(W_2)$ are 
\eqna
F_1 &=& x_2' - \DS{\frac{1}{2}x_{1,1}' +\frac{1}{8} x_{\ubar,\ubar}'}\\
F_2&=& \DS{\frac{1}{2} (x_{\dbar}' + x_{1,\ubar}')}\\
F_3 &=& \DS{\frac{1}{2} x_{1,1}'} \\
F_4 &=& \DS{\frac{1}{8} x_{\ubar,\ubar}'}
\endeqna

If we denote by $A_i$ the block $\QM\Sigma'(W_2) F_i$, then 
$A_i \simeq \QM$ if $i \in \{1,3,4\}$ and 
$A_2$ is isomorphic to the algebra to upper triangular $2 \times 2$-matrices 
(see \cite[\S 6]{BH}). In particular, $\QM\Sigma'(W_2)$ is hereditary. 

For information, we provide the dimensions of the left ideal, right ideal, 
two-sided ideal generated by $x_C$ (for $C \in \Comp(2)$) and also the dimension 
of the centralizer of $x_C$. In this 
table, $A$ denotes the algebra $\QM\Sigma'(W_2)$. 

$$\begin{array}{|c||c|c|c|c|c|c|}
\hline
\vertical  x & x_2 & x_{\dbar} & x_{1,1} & x_{1,\ubar} & x_{\ubar,1} & x_{\ubar,\ubar}   \\
\hline
\hline
\vertical \dim_\QM Ax     & 6 & 3 & 3 & 2 & 2 & 1 \\
\hline
\vertical \dim_\QM xA     & 6 & 2 & 4 & 3 & 3 & 1 \\
\hline
\vertical \dim_\QM AxA    & 6 & 3 & 4 & 3 & 3 & 1 \\
\hline
\vertical \dim_\QM Z_A(x) & 6 & 5 & 5 & 5 & 5 & 6 \\
\hline
\end{array}$$

\bigskip

\soussection{The case ${\boldsymbol{n=3}}$} 
The character table of $\QM\Sigma'(W_3)$ is:

$$\begin{array}{|c|cccccccccc|}
\hline
\vertical & x_3 & x_{\tbar} & x_{2,1} & x_{2,\ubar} & x_{1,\dbar} & x_{\dbar,\ubar} & 
x_{1,1,1} & x_{1,1,\ubar} & x_{1,\ubar,\ubar} & x_{\ubar,\ubar,\ubar} \\
\hline
\vertical \pi_{3;\vide}^\QM    & 1 & . & . & . & . & . & . & . & . & .\\
\vertical \pi_{\vide;3}^\QM    & 1 & 2 & . & . & . & . & . & . & . & .\\
\vertical \pi_{21;\vide}^\QM   & 1 & . & 1 & . & . & . & . & . & . & .\\
\vertical \pi_{2;1}^\QM        & 1 & . & 1 & 2 & . & . & . & . & . & .\\
\vertical \pi_{1;2}^\QM        & 1 & . & 1 & . & 2 & . & . & . & . & .\\
\vertical \pi_{\vide;21}^\QM   & 1 & 4 & 1 & 2 & 2 & 4 & . & . & . & .\\
\vertical \pi_{111;\vide}^\QM  & 1 & . & 3 & . & . & . & 6 & . & . & .\\
\vertical \pi_{11;1}^\QM       & 1 & . & 3 & 2 & . & . & 6 & 4 & . & .\\
\vertical \pi_{1;11}^\QM       & 1 & . & 3 & 4 & 4 & . & 6 & 8 & 8 & .\\
\vertical \pi_{\vide;111}^\QM  & 1 & 8 & 3 & 6 & 12 & 24 & 6 & 12 & 24 & 48\\ 
\hline
\end{array}$$

\bigskip

We can take for the family $(E_\l^\QM)_{\l \in \Bip(3)}$ the following 
idempotents:
\eqna
E_{3;\vide}^\QM &=& 
\DS{x_3'- x_{1,2}' + \frac{1}{4}x_{\ubar,\dbar}' + \frac{1}{3}x_{1,1,1}'
- \frac{1}{6} x_{\ubar,1,\ubar}' +\frac{1}{12}(x_{1,\ubar,\ubar}'+x_{\ubar,\ubar,1}')} \\
E_{\vide;3}^\QM &=& \DS{\frac{1}{2}(x_{\tbar}' + x_{2,\ubar}'-x_{\ubar,2}') 
-\frac{1}{3}x_{1,1,\ubar}'+\frac{1}{6} (x_{1,\ubar,1}' + x_{\ubar,1,1,}')}\\
E_{21;\vide}^\QM &=& \DS{x_{1,2}' - \frac{1}{2}x_{1,1,1}' + \frac{1}{8}x_{1,\ubar,\ubar}'}\\
E_{2;1}^\QM&=& \DS{\frac{1}{2} x_{\ubar,2}' -\frac{1}{4}x_{\ubar,1,1}'+\frac{1}{16}
x_{\ubar,\ubar,\ubar}'}\\
E_{1;2}^\QM &=& \DS{\frac{1}{2} x_{1,\dbar}' +\frac{1}{4}(x_{1,1,\ubar}'- x_{1,\ubar,1}')} \\
E_{\vide;21}^\QM&=& \DS{\frac{1}{4} x_{\dbar,\ubar}' + 
\frac{1}{8}(x_{1,\ubar,1}'-x_{\ubar,1,\ubar}')}\\
E_{111;\vide}^\QM&=&\DS{\frac{1}{6} x_{1,1,1}'}\\
E_{11;1}^\QM&=& \DS{\frac{1}{12}(x_{1,1,\ubar}' + x_{1,\ubar,1}' + x_{\ubar,1,1}')}\\
E_{1;11}^\QM &=& \DS{-\frac{1}{12} x_{1,\ubar,\ubar}' + \frac{7}{24} x_{\ubar,1,\ubar}' -
\frac{1}{12} x_{\ubar,\ubar,1}' } \\
E_{\vide;111}^\QM&=& \DS{\frac{1}{48} x_{\ubar,\ubar,\ubar}'}\\
\endeqna

\bigskip

The Cartan matrix of $\QM\Sigma'(W_3)$ is given by

\bigskip

$$\begin{array}{|c||cccc|cccc|c|c|}
\hline
\vertical & \DC_{3;\vide}^\QM & \DC_{21;\vide}^\QM & \DC_{\vide;21}^\QM & \DC_{1;11}^\QM & 
\DC_{\vide;3}^\QM & \DC_{2;1}^\QM & \DC_{1;2}^\QM & \DC_{11;1}^\QM & \DC_{111;\vide}^\QM & 
\DC_{\vide;111}^\QM\\
\hline
\hline
\vertical \PC_{3;\vide}^\QM   & 1 & 1 & 1 & 1 & . & . & . & . & . & .\\
\vertical \PC_{21;\vide}^\QM  & . & 1 & . & . & . & . & . & . & . & .\\
\vertical \PC_{\vide;21}^\QM  & . & . & 1 & 1 & . & . & . & . & . & .\\
\vertical \PC_{1;11}^\QM      & . & . & . & 1 & . & . & . & . & . & .\\
\hline
\vertical \PC_{\vide;3}^\QM   & . & . & . & . & 1 & 1 & 1 & 1 & . & .\\
\vertical \PC_{2;1}^\QM       & . & . & . & . & . & 1 & . & . & . & .\\
\vertical \PC_{1;2}^\QM       & . & . & . & . & . & . & 1 & 1 & . & .\\
\vertical \PC_{11;1}^\QM      & . & . & . & . & . & . & . & 1 & . & .\\
\hline
\vertical \PC_{111;\vide}^\QM & . & . & . & . & . & . & . & . & 1 & .\\
\hline
\vertical \PC_{\vide;111}^\QM & . & . & . & . & . & . & . & . & . & 1\\ 
\hline
\end{array}$$

\bigskip

The primitive central idempotents of $\QM\Sigma'(W_3)$ are 
\eqna
F_1 &=& \DS{x_3' +\frac{1}{4}(x_{\dbar,\ubar}'+ x_{\ubar,\dbar}'+x_{1,\ubar,\ubar}') 
- \frac{1}{6} x_{1,1,1}'}\\
F_2&=& \DS{\frac{1}{2}( x_{\tbar}' + x_{2,\ubar}' + x_{1,\dbar}') + 
\frac{5}{48} x_{\ubar,\ubar,\ubar}'}\\
F_3 &=& \DS{\frac{1}{6} x_{1,1,1}'}\\
F_4&=&\DS{\frac{1}{48} x_{\ubar,\ubar,\ubar}'}
\endeqna

For information, we provide the dimensions of the left ideal, right ideal, 
two-sided ideal generated by $x_C$ (for $C \in \Comp(3)$) and also the dimension 
of the centralizer of $x_C$. In these 
tables, $A$ denotes the algebra $\QM\Sigma'(W_3)$.

$$\begin{array}{|c||c|c|c|c|c|c|c|c|c|c|}
\hline
\vertical  x & x_3 & x_{\tbar} & x_{2,1} & x_{1,2} & x_{2,\ubar} & x_{\ubar,2} & 
x_{\dbar,1} & x_{1,\dbar} & x_{\dbar,\ubar} & x_{\ubar,\dbar} \\
\hline
\hline
\vertical \dim_\QM Ax     &18& 7&10&10& 6& 6& 6& 6& 3& 3 \\
\hline
\vertical \dim_\QM xA     &18& 4&16&16&11&11& 8& 8& 3& 3 \\
\hline
\vertical \dim_\QM AxA    &18& 9&16&16&11&11&10&10& 5& 5 \\
\hline
\vertical \dim_\QM Z_A(x) &18&13&10&10&12&12&13&13&16&16 \\
\hline
\end{array}$$

$$\begin{array}{|c||c|c|c|c|c|c|c|c|}
\hline
\vertical  x & x_{1,1,1} & 
x_{1,1,\ubar} & x_{1,\ubar,1} & x_{\ubar,1,1} & x_{1,\ubar,\ubar} & x_{\ubar,1,\ubar} 
& x_{\ubar,\ubar,1} & x_{\ubar,\ubar,\ubar} \\
\hline
\hline
\vertical \dim_\QM Ax     & 4& 3& 3& 3& 2& 2& 2& 1 \\
\hline
\vertical \dim_\QM xA     & 8& 7& 7& 7& 4& 4& 4& 1 \\
\hline
\vertical \dim_\QM AxA    & 8& 7& 7& 7& 4& 4& 4& 1 \\
\hline
\vertical \dim_\QM Z_A(x) &14&14&14&14&16&16&16&18 \\
\hline
\end{array}$$

\bigskip

\soussection{The case ${\boldsymbol{n=4}}$} 
We shall give here only the Cartan matrix and the central 
idempotents of $A$. The Cartan matrix is given by

\bigskip

{\small
\begin{centerline}{$\begin{array}{|c||cccccccc|ccccccccc|c|c|c|}
\hline
&4     &31   &\vide&\vide&211  &2 &1 &11&\vide&3&1&2&21&11&\vide&111&1  &22   &1111 &\vide\\
&\vide &\vide&31   &22   &\vide&11&21&11&4    &1&3&2&1 &2 &211  &1  &111&\vide&\vide&1111\\
\hline 
\hline
\vertical \PC_{\!4;\vide}^\QM  
&~1~&~1~&~1~&~.~&~1~&~1~&~2~&~1~&~.~&~.~&~.~&~.~&~.~&~.~&~.~&~.~&~.~&~.~&~.~&~.~\\
\vertical \PC_{\!31;\vide}^\QM             &.&1&.&.&1&.&1&1&.&.&.&.&.&.&.&.&.&.&.&.\\
\vertical \PC_{\!\vide;31}^\QM             &.&.&1&.&.&1&1&1&.&.&.&.&.&.&.&.&.&.&.&.\\
\vertical \PC_{\!\vide;22}^\QM             &.&.&.&1&.&.&1&1&.&.&.&.&.&.&.&.&.&.&.&.\\
\vertical \PC_{\!211;\vide}^\QM            &.&.&.&.&1&.&.&.&.&.&.&.&.&.&.&.&.&.&.&.\\
\vertical \PC_{\!2;11}^\QM                 &.&.&.&.&.&1&.&.&.&.&.&.&.&.&.&.&.&.&.&.\\
\vertical \PC_{\!1;21}^\QM                 &.&.&.&.&.&.&1&1&.&.&.&.&.&.&.&.&.&.&.&.\\
\vertical \PC_{\!1\!1;1\!1}^\QM            &.&.&.&.&.&.&.&1&.&.&.&.&.&.&.&.&.&.&.&.\\
\hline
\vertical \PC_{\!\vide;4}^\QM              &.&.&.&.&.&.&.&.&1&1&1&1&2&1&1&1&1&.&.&.\\
\vertical \PC_{\!3;1}^\QM                  &.&.&.&.&.&.&.&.&.&1&.&.&1&.&1&.&1&.&.&.\\
\vertical \PC_{\!1;3}^\QM                  &.&.&.&.&.&.&.&.&.&.&1&.&1&1&.&1&.&.&.&.\\
\vertical \PC_{\!2;2}^\QM                  &.&.&.&.&.&.&.&.&.&.&.&1&1&.&.&.&.&.&.&.\\
\vertical \PC_{\!21;1}^\QM                 &.&.&.&.&.&.&.&.&.&.&.&.&1&.&.&.&.&.&.&.\\
\vertical \PC_{\!1\!1;2}^\QM               &.&.&.&.&.&.&.&.&.&.&.&.&.&1&.&1&.&.&.&.\\
\vertical \PC_{\!\vide;21\!1}^\QM          &.&.&.&.&.&.&.&.&.&.&.&.&.&.&1&.&1&.&.&.\\
\vertical \PC_{\!1\!1\!1;1}^\QM            &.&.&.&.&.&.&.&.&.&.&.&.&.&.&.&1&.&.&.&.\\
\vertical \PC_{\!1;1\!1\!1}^\QM            &.&.&.&.&.&.&.&.&.&.&.&.&.&.&.&.&1&.&.&.\\
\hline
\vertical \PC_{\!22;\vide}^\QM             &.&.&.&.&.&.&.&.&.&.&.&.&.&.&.&.&.&1&.&.\\
\hline
\vertical \!\PC_{\!1\!1\!1\!1;\vide}^\QM\! &.&.&.&.&.&.&.&.&.&.&.&.&.&.&.&.&.&.&1&.\\
\hline
\vertical \!\PC_{\!\vide;1\!1\!1\!1}^\QM\! &.&.&.&.&.&.&.&.&.&.&.&.&.&.&.&.&.&.&.&1\\ 
\hline
\end{array}$}\end{centerline}}

\bigskip

The primitive central idempotents of $\QM\Sigma'(W_4)$ are
\eqna
F_1&=& \DS{x_4' -\frac{1}{2} x_{2,2}' + \frac{1}{4}(x_{\tbar,\ubar}+x_{\ubar,\tbar}' 
+x_{2,1,1}'+x_{1,1,2}'+x_{\dbar,\dbar}'+x_{1,\ubar,\dbar}'+x_{1,\dbar,\ubar}')}\\
&& \DS{+ \frac{1}{6} x_{1,1,1,1}' +\frac{3}{16} x_{2,\ubar,\ubar}' 
- \frac{1}{16} x_{\ubar,\ubar,2}'+\frac{1}{32}(x_{1,1,\ubar,\ubar}'+x_{\ubar,\ubar,1,1}') 
+ \frac{5}{96}x_{\ubar,\ubar,\ubar,\ubar}'}\\
F_2&=&\DS{\frac{1}{2}(x_{\bar{4}}'+x_{3,\ubar}'+x_{1,\tbar}'+x_{2,\dbar}') 
+ \frac{1}{8}(x_{\dbar,\ubar,\ubar}'+x_{\ubar,\dbar,\ubar}'+x_{\ubar,\ubar,\dbar}'+
x_{1,\ubar,\ubar,\ubar}')}\\
F_3&=&\DS{\frac{1}{2} x_{2,2}' -\frac{1}{4}(x_{2,1,1}'+x_{1,1,2}') + \frac{1}{8} x_{1,1,1,1}' 
+\frac{1}{16}(x_{2,\ubar,\ubar}'+x_{\ubar,\ubar,2}')}\\
&& \DS{-\frac{1}{32} (x_{1,1,\ubar,\ubar}'+x_{\ubar,\ubar,1,1}') 
+ \frac{1}{128} x_{\ubar,\ubar,\ubar,\ubar}'}\\
F_4&=& \DS{\frac{1}{24} x_{1,1,1,1}'}\\
F_5&=& \DS{\frac{1}{384} x_{\ubar,\ubar,\ubar,\ubar}'}
\endeqna

\section{Questions}

\medskip

Let us raise here some questions about the representation 
theory of the Mantaci-Reutenauer algebra $K\Sigma'(W_n)$:

\medskip

(1) Determine the centre of 
$K\Sigma'(W_n)$, or at least its dimension (in characteristic 
zero, its dimension determines its structure because it is 
split semisimple by Corollary \ref{centre semisimple}).

\medskip

(2) Compute the Cartan matrix of $\QM\Sigma'(W_n)$. Note that the Theorem 
\ref{cartan induction} provides a first induction argument. 

\medskip

($2^+$) Determine the Loewy series of the projective indecomposable $K\Sigma'(W_n)$-modules. 
Determine the Loewy length of 
$\FM_{\! 2}\Sigma'(W_n)$ (see Remark \ref{l2}: it is probably equal 
to $2n-1$ if $n \ge 2$). 

\medskip

(3) For which values of $n$ is the algebra $\QM\Sigma'(W_n)$ hereditary? 
It is reasonable to expect that it is hereditary if and only if $n \in \{1,2,3\}$. 
Note that it is not hereditary for $n \in \{4,5\}$. 

\medskip

($3^+$) Compute the path algebra of $\QM\Sigma'(W_n)$.

\medskip

(4) Is the inclusion (4) in the proof of Proposition \ref{loewy groth} always an equality? 
For the analogous statement for the symmetric group, we have 
an equality \cite[Theorem A]{bonnafe sym}.

\medskip

(5) Is the morphism $\Res_k^n : \ZM\Sigma'(W_n) \to \ZM\Sigma'(W_k)$ surjective? 
Compare with Proposition \ref{surjectif B}. 
Note that the morphism $\Res_{W_k}^{W_n} : \ZM\Irr W_n \to \ZM\Irr W_k$ 
is surjective.

\medskip

(6) The Corollary \ref{dimension proj} suggests, by analogy with the case 
of the symmetric group, the following question:
if $\l \in \Bip(n)$, does there exist a 
linear character $\z_\l$ of $C_{W_n}(\cox_\l)$ such that 
$\CM W_n E_\l^\CM$ affords the character $\Ind_{C_{W_n}(\cox_\l)}^{W_n} \z_\l$? 
In fact, the answer to this question is negative in general, even for $n=2$ 
(take $\l=((2);\vide)$). Computations using {\tt CHEVIE} (for $n \le 4$) 
suggests that the following slight modification of the previous question 
could have a positive answer: 
if $\l \in \Bip(n)$, does there exist a 
linear character $\z_\l$ of $C_{W_n}(\cox_{\l^{\mathrm{opp}}})$ such that 
$\CM W_n E_\l^\CM$ affords the character 
$\Ind_{C_{W_n}(\cox_{\l^{\mathrm{opp}}})}^{W_n} \z_\l$? 
Here, if $\l=(\l^+,\l^-)$, we have denoted by $\l^{\mathrm{opp}}$ the 
bipartition $(\l^-,\l^+)$. 

\bigskip

\noindent{\sc Remark - } 
It is readily seen that $|C_{W_n}(w_\l)|=|C_{W_n}(w_{\l^{\mathrm{opp}}})|$.\finl

\bigskip

\example{exemple question 6} 
We keep here the notation of Example \ref{x un un prime}. 
Let $r \in \{0,1,2,\dots,n\}$. Then 
$$\Cb(I_r^+)=(\underbrace{1,1,\dots,1}_{\text{$r$ times}}, 
\underbrace{-1,-1,\dots,-1}_{\text{$n-r$ times}})$$
and, if we set $\l(r)=\lamb(\Cb(I_r^+))$, then 
$$\l(r)=(\underbrace{(1,1,\dots,1)}_{\text{$r$ times}}, 
\underbrace{(1,1,\dots,1)}_{\text{$n-r$ times}}).$$
We shall prove that the answer to question 6 is positive 
whenever $\l=\l(r)$. 

Let us make this statement more precise. 
We may choose for $\cox_{\l(r)^{\mathrm{opp}}}$ 
the element $\cox_{\Cb(I_r^+)}=t_{r+1}t_{r+2}\dots t_n$. 
Then 
$$C_{W_n}(\cox_{\l(r)^{\mathrm{opp}}}) = W_{r,n-r}.$$
Let $\g_r \boxtimes 1_{n-r}$ denote the linear character of 
$W_{r,n-r} \simeq W_r \times W_{n-r}$ which is equal to 
$\g_r$ on the component $W_r$ and which is trivial on $W_{n-r}$ 
(recall that the linear character $\g_r$ of $W_r$ has been defined 
in Example \ref{x un un prime}). Then: 
\equat\label{question 6}
\text{\it The $\QM W_n$-module $\QM W_n E_{\l(r)}$ affords the character 
$\Ind_{W_{r,n-r}}^{W_n} (\g_r \boxtimes 1_{n-r})$.}
\endequat
Let us prove \ref{question 6}. First, note that $\g_r \boxtimes 1_{n-r}$ 
is just the restriction of the map $\g_{I_r^+} : W_n \to \{1,-1\}$ 
to $W_{r,n-r}$ defined in Example \ref{x un un prime}. 
By Example \ref{x un un prime idempotent}, we may take 
for $E_{\l(r)}^\QM$ the idempotent $x_{I_r^+}'/(2^{n-r} r! (n-r)!)$. 
Let $e_{\SG_n}=(1/n!)\sum_{\s \in \SG_n} \s$ and let 
$$E(r)=\frac{1}{|W_{r,n-r}|} \sum_{\s \in W_{r,n-r}} (\g_r \boxtimes 1_{n-r})(\s) \s.$$
Then the module $\QM W_n E(r)$ affords the character 
$\Ind_{W_{r,n-r}}^{W_n} (\g_r \boxtimes 1_{n-r})$ and, by Example 
\ref{x un un prime idempotent}, we have 
$E_{\l(r)} = (|W_n|/|W_{r,n-r}) e_{\SG_n} E(r)$. 
Therefore, $\QM W_n E_{\l(r)} \subset \QM W_n E(r)$. But 
$\dim_\QM \QM W_n E(r)=|W_n|/|W_{r,n-r}|$ and $\dim_\QM \QM W_n E_{\l(r)} 
= |\CC(\l(r))|=|W_n|/|W_{r,n-r}|$, so we get that 
$\QM W_n E_{\l(r)} = \QM W_n E(r)$.\finl

\end{document}